\documentclass[12pt]{article}

\usepackage{psfig}
\usepackage{amsfonts}
\newtheorem{lemma}{Lemma}[section]
\newtheorem{theorem}{Theorem}[section]
\newtheorem{conjecture}{Conjecture}[section]

\begin{document}

\begin{center}
{\bf \Large New Conjectural Lower Bounds on the Optimal Density of Sphere Packings}
\end{center}
\bigskip

\begin{center}
S. Torquato$^{1,2,3,4}$ and F. H. Stillinger$^1$
 
$^{1}$Department of Chemistry, $^2$Program in Applied and Computational Mathematics,
$^3$Princeton Institute for
 the Science and Technology of Materials, and \\
$^4$ Princeton Center for Theoretical Physics\\ 
Princeton University, Princeton, NJ 08544 \\
\end{center}

\centerline{\bf ABSTRACT}

Sphere packings in high dimensions interest  mathematicians and
physicists and have direct applications in communications theory.
Remarkably, no one has been able to provide exponential improvement on a
100-year-old lower bound on the maximal packing density due to Minkowski in 
$d$-dimensional Euclidean space $\mathbb{R}^d$. The asymptotic behavior of this bound 
is controlled by $2^{-d}$ in high dimensions.  Using an optimization procedure 
that we introduced earlier \cite{To02c} and a conjecture
concerning the existence of disordered sphere packings in $\mathbb{R}^d$, we obtain
a conjectural lower bound on the density whose asymptotic behavior
is controlled by $2^{-0.77865\ldots d}$, thus providing the putative
exponential improvement of Minkowski's bound. The conjecture states that a hard-core
nonnegative tempered distribution is a pair correlation function of a
translationally invariant disordered sphere packing in $\mathbb{R}^d$ for
asymptotically large $d$ if and only if the Fourier transform of the
autocovariance function is nonnegative. The conjecture is supported by two
explicit analytically characterized disordered packings, numerical
packing constructions in low dimensions, known necessary conditions that only
have relevance in very low dimensions, and the fact that we can
recover the forms of known rigorous lower bounds. A byproduct of our approach is an 
asymptotic conjectural lower bound on the average kissing number whose behavior is 
controlled by $2^{0.22134\ldots d}$, which is to be compared to the best known 
asymptotic lower bound on the individual kissing number of $2^{0.2075\ldots d}$.
Interestingly, our optimization procedure is precisely the dual of a primal linear program devised by Cohn and 
Elkies \cite{Co03} to obtain upper bounds on the density, and hence
has implications for linear programming bounds. This connection
proves that our density estimate can never exceed the Cohn-Elkies upper
bound, regardless of the validity of our conjecture.

\newpage
\section{Introduction}

A collection of congruent spheres in $d$-dimensional Euclidean space $\mathbb{R}^d$ 
is called a sphere packing $P$ 
if no two of the spheres have an interior point in common.
The {\it packing density} or simply density $\phi(P)$ of a sphere packing is the fraction of
space $\mathbb{R}^d$ covered by the spheres. We will call
\begin{equation}
\phi_{\mbox{\scriptsize max}}= \sup_{P\subset \mathbb{R}^d} \phi(P)
\end{equation}
the {\it maximal density}, where the supremum is taken over
all packings in $\mathbb{R}^d$.
The sphere packing problem seeks to answer the following 
question: Among all packings of congruent spheres,
what is the maximal packing density $\phi_{\mbox{\scriptsize max}}$, i.e., largest
fraction of $\mathbb{R}^d$ covered by the spheres,
and what are the corresponding arrangements of the spheres \cite{Ro64,Co93}?
The sphere packing problem is of great fundamental and practical interest,
and arises in a variety of contexts, including classical ground states
of  matter in low dimensions \cite{Chaik95}, 
the famous Kepler conjecture for $d=3$ \cite{Ha98}, error-correcting
codes \cite{Co93} and spherical codes \cite{Co93}.

For arbitrary $d$, the sphere packing problem  is notoriously difficult to solve.
In the case of packings of congruent $d$-dimensional
spheres, the exact solution is known for the first
three space dimensions. For $d=1$, the answer is trivial
because the spheres tile the space so that $\phi_{\mbox{\scriptsize max}}=1$.
In two dimensions, the optimal solution is the triangular lattice 
arrangement (also called the hexagonal packing) with $\phi_{\mbox{\scriptsize max}}=\pi/\sqrt{12}$.
In three dimensions, the Kepler conjecture
that the face-centered cubic lattice arrangement
provides the densest packing with
$\phi_{\mbox{\scriptsize max}}=\pi/\sqrt{18}$ 
was only recently proved by Hales \cite{Ha98}.
For $3< d <10$, 
the densest known packings of congruent spheres are lattice
packings (defined below). However, for sufficiently large $d$,
lattice packings are likely not to be the densest. Indeed, this paper suggests
that disordered sphere arrangements may be the
densest packings as $d \rightarrow \infty$.

We review some basic definitions that we will subsequently employ.
A packing is {\it saturated} if there is no space available to add another sphere
 without overlapping the existing particles. 
The set of lattice packings is a subset of the set of sphere
packings in $\mathbb{R}^d$. A {\it lattice} $\Lambda$ in $\mathbb{R}^d$ is a subgroup
consisting of the integer linear combinations of vectors that constitute a basis for $\mathbb{R}^d$.
A {\it lattice packing} $P_L$ is one in which  the centers of nonoverlapping spheres
are located at the points of $\Lambda$. 
In a lattice packing, the space $\mathbb{R}^d$ can be geometrically divided into identical
regions $F$ called {\it fundamental cells}, each of which contains the center
of just one sphere. Thus, the density of a lattice packing $\phi^L$ consisting
of spheres of unit diameter is given by
\begin{equation}
\phi^L= \frac{v_1(1/2)}{\mbox{Vol}(F)},
\end{equation}
where
\begin{equation}
v_1(R) = \frac{\pi^{d/2}}{\Gamma(1+d/2)} R^d
\label{v(R)}
\end{equation}
is the volume of a $d$-dimensional sphere of radius $R$ and
$\mbox{Vol}(F)$ is the volume of a fundamental cell.
We will call 
\begin{equation}
\phi^L_{\mbox{\scriptsize max}}= \sup_{P_L\subset \mathbb{R}^d} \phi(P_L)
\end{equation}
the maximal density among all lattice packings in $\mathbb{R}^d$.
For a general packing of spheres of unit diameter for which a density $\phi(P)$ exists, it is
useful to introduce the {\it number} (or {\it center}) density $\rho$
defined by
\begin{equation}
\rho = \frac{\phi(P)}{v_1(1/2)},
\end{equation}
which therefore can be interpreted as the average number of sphere centers
per unit volume.

%A more general notion than a lattice packing is a periodic
%packing. A {\it periodic} packing of congruent spheres
%is obtained by placing a fixed nonoverlapping configuration of $N$ particles (where $N\ge 1$)
%in each fundamental cell of a lattice $\GL$. Thus, the packing is still
%periodic under translations by $\GL$, but the $N$ spheres can occur
%anywhere in the chosen fundamental cell subject to the nonoverlap condition.
%The packing density of a such periodic packing in which the spheres 
%have unit diameter is given by
%\begin{equation}
%\phi=\frac{N v_1(1/2)}{V}=\rho v_1(1/2),
%\end{equation}
%where 
%\begin{equation}
%v_1(R) = \frac{\pi^{d/2}}{\Gamma(1+d/2)} R^d 
%\label{v(R)}
%\end{equation}
%is the volume of a $d$-dimensional sphere of radius $R$ and $V$ is the volume
%of the fundamental cell. Here the number (or center) density $\rho=N/V$ 
%is just the number of spheres per unit volume.

%In some cases, the configuration of spheres within the fundamental
%cell of a periodic packing may be described by a probability measure.
%We will be particularly interested in such disordered arrangements in the limit
%that $N \rightarrow \infty$ and $V \rightarrow \infty$ such
%that the number density $\rho$ exists. We
%call such an infinite packing a disordered sphere packing.

Three distinct categories of packings have been
distinguished, depending on their behavior with respect to nonoverlapping
geometric constraints and/or externally imposed virtual displacements:
{\it locally} jammed, {\it collectively} jammed, and {\it strictly} jammed \cite{To01b,To03c,Do04a}.
Loosely speaking, these {\it jamming} categories, listed in order of increasing stringency,
reflect the degree of mechanical stability of the packing.
A packing is locally jammed if each particle in the system is locally trapped by
its neighbors; i.e., it cannot be translated while fixing the positions
of all other particles. Each sphere simply has to have at least $d+1$
contacts with neighboring spheres, not all in the same $d$-dimensional
hemisphere. A collectively jammed packing is any locally jammed configuration in which no
finite subset of particles can simultaneously be continuously displaced so that its members
move out of contact with one another and with the remainder set.
A strictly jammed packing is any collectively jammed configuration that disallows
all globally uniform volume-nonincreasing deformations of the system
boundary. Importantly, the jamming category is generally dependent
on the type of boundary conditions imposed (e.g., hard wall or periodic
boundary conditions) as well as the shape of the boundary.  The range of possible densities 
for a given jamming category decreases with increasing stringency of the category. 
Whereas the lowest-density states of collectively and strictly jammed packings
in two or three dimensions are currently unknown, one can achieve
locally jammed packings with vanishing density \cite{Bor64}.
This classification of packings according to jamming criteria is closely related 
to the concepts of {}``rigid'' and ``stable'' packings found in the mathematics literature
\cite{Co98}.

In the next section, we summarize some previous upper and lower bounds
on the maximal density. For large $d$, the asymptotic behavior of the well-known Minkowski 
lower bound \cite{Mi05} on the maximal density is controlled by $2^{-d}$.
Thus far, no one has been able to provide exponential improvement on
this lower bound. Using an optimization procedure and a conjecture
concerning the existence of disordered sphere packings in high dimensions,
we obtain conjectural lower bounds that yield the long-sought 
asymptotic exponential improvement on Minkowski's bound. We believe that consideration 
of truly disordered packings is the key notion that will yield exponential improvement
on Minkowski's lower bound. A byproduct of our approach is an 
asymptotic conjectural lower bound on the average kissing number that is superior           
to the best known asymptotic lower bound on the individual kissing number.

\section{Some Previous Upper and Lower Bounds}

The {\it nonconstructive} lower bounds of Minkowski \cite{Mi05} established the
existence of reasonably dense lattice packings. He found
that the maximal density $\phi^L_{\mbox{\scriptsize max}}$ among all lattice packings 
for $d \ge 2$ satisfies
\begin{equation}
\phi^L_{\mbox{\scriptsize max}} \ge \frac{\zeta(d)}{2^{d-1}},
\label{mink}
\end{equation}
where $\zeta(d)=\sum_{k=1}^\infty k^{-d}$ is the Riemann zeta function.
Note that for large values of $d$,
the asymptotic behavior of the Minkowski lower bound is controlled by $2^{-d}$.
Since 1905, many extensions and generalizations of (\ref{mink})
have been obtained \cite{Da47,Ball92,Co93}, but none of these investigations have been able to improve
upon the dominant exponential term $2^{-d}$. It is useful to note that
the density of a saturated packing of congruent spheres
in $\mathbb{R}^d$ for all $d$ satisfies 
\begin{equation}
\phi \ge \frac{1}{2^d}.
\label{sat}
\end{equation}
The proof is trivial. A saturated packing of congruent spheres
of unit diameter  and density $\phi$ in $\mathbb{R}^d$ has the property that each point in space lies
within a unit distance from the center of some sphere. Thus, a covering
of the space is achieved if each sphere center is encompassed by a sphere
of unit radius  and the density of this covering is $2^d \phi \ge 1$.
Thus, the bound (\ref{sat}), which is sometimes called the ``greedy"
lower bound, has same  the dominant exponential term as (\ref{mink}).
In Section 4.1, we show that there exists a construction of
a disordered packing of congruent spheres that realizes the weaker 
lower bound  of (\ref{sat}), i.e., $\phi=2^{-d}$. 

The best currently known lower
bound on $\phi^L_{\mbox{\scriptsize max}}$ was obtained by Ball \cite{Ball92}.
He found that 
\begin{equation}
\phi^L_{\mbox{\scriptsize max}} \ge \frac{2(d-1)\zeta(d)}{2^{d}}.
\label{ball}
\end{equation}
As $d \rightarrow \infty$, 
Table \ref{table-lower} gives the dominant asymptotic behavior of several
lower bounds on $\phi^L_{\mbox{\scriptsize max}}$ for  large $d$.

\begin{table}[bthp]
\centering
\caption{ Dominant asymptotic behavior of lower bounds on $\phi_{\mbox{\scriptsize max}}^L$
for  large $d$.\label{table-lower}}
\begin{tabular}{c|c}
\multicolumn{2}{c}{~}\\ \hline\hline
$(2)  2^{-d}$  &  \mbox{Minkowski (1905)} \\ \hline

$[\ln (\sqrt{2}) d] 2^{-d}$  &   \mbox{Davenport and Rogers (1947)}  \\ \hline

$(2 d)2^{-d}$  &   \mbox{Ball (1992)}   \\ \hline\hline
\end{tabular}
\end{table}

Nontrivial upper bounds on the
maximal density $\phi_{\mbox{\scriptsize max}}$ for any sphere
packing in $\mathbb{R}^d$ have been derived. 
Blichfeldt \cite{Bl29} showed that the maximal packing density for all $d$ satisfies
$\phi_{\mbox{\scriptsize max}} \le (d/2+1)2^{-d/2}.$
This upper bound was improved by Rogers \cite{Ro58,Ro64}
by an analysis of the Voronoi cells. For large $d$, Rogers'
upper bound asymptotically becomes $d 2^{-d/2}/e$. Kabatiansky and Levenshtein
\cite{Ka78} found an even stronger bound, which in the limit $d\rightarrow \infty$
yields $\phi_{\mbox{\scriptsize max}} \le 2^{-0.5990d(1+o(1))}$.
Cohn and Elkies \cite{Co03} obtained and computed linear programming upper bounds,
which provided improvement over Rogers' upper bound for dimensions
4 through 36. They also conjectured that their approach could be used
to prove sharp bounds in 8 and 24 dimensions. 
Indeed, Cohn and Kumar \cite{Co04} used these techniques
to prove that the Leech lattice is the unique densest lattice in $\Re^{24}$. They
also proved that no sphere packing in
$\Re^{24}$ can exceed the density of the Leech lattice
by a factor of more than $1+1.65 \times 10^{-30}$, 
and gave a new proof that $E_8$ is
the unique densest lattice in $\Re^8$.
Table \ref{upper} provides the dominant asymptotic behavior of several
upper bounds on $\phi_{\mbox{\scriptsize max}}$ for  large $d$.

\begin{table}[bthp]
\centering
\caption{ Dominant asymptotic behavior of upper bounds on $\phi_{\mbox{\scriptsize max}}$
for large $d$.\label{upper}}
\begin{tabular}{c|c}
\multicolumn{2}{c}{~}\\ \hline\hline
$(d/2)2^{-0.5d}$  &  \mbox{Blichfeldt (1929)} \\ \hline

  $(d/e)2^{-0.5d}$  &   \mbox{Rogers (1958)}  \\ \hline

  $2^{-0.5990d}$  &   \mbox{Kabatiansky and Levenshtein (1979)}   \\ \hline\hline
\end{tabular}
\end{table}

\section{Realizability of Point Processes}
\label{realize}

As will be described in Section \ref{new}, our new approach to lower bounds on the density of
sphere packings in $\mathbb{R}^d$ rests on whether certain one- and two-point correlation
functions are realizable by sphere packings. As we will discuss, a sphere packing
can be regarded to be a special case of a point process and so
a more general question concerns the necessary and sufficient conditions
for the realizability of point processes in $\mathbb{R}^d$. Before discussing the realizability 
of  point processes, it is useful to recall some basic results from the theory of
spatial stochastic (or random) processes.  Let ${\bf x} \equiv (x_1,x_2,\ldots x_d)$ 
denote a vector position in $\mathbb{R}^d$. Consider a stochastic process 
$\{Y({\bf x};\omega) : {\bf x}\in \mathbb{R}^d; \omega\in \Omega\}$, where  
$Y({\bf x};\omega)$ is a real-valued random variable, 
$\omega$ is a realization generated by the stochastic
process and ($\Omega,{\cal F},{\cal P}$) is a probability space
(i.e., $\Omega$ is a sample space, $\cal F$ is a $\sigma$-algebra
of measurable subsets of $\Omega$, and $\cal P$ is a probability measure).
For simplicity, we will often suppress the variable $\omega$.

\subsection{Ordinary Stochastic Processes}

We will assume that the stochastic process
is translationally invariant (i.e., statistically
homogeneous in space). Let us further assume that the mean 
$\mu=\langle Y({\bf x}) \rangle$  and autocovariance function 
\begin{equation}
\chi({\bf r})=\langle [Y({\bf x})-\mu][Y({\bf x}+{\bf r})-\mu]\rangle
\label{auto}
\end{equation}
exist, where angular brackets denote an expectation, i.e., an average
over all realizations. The fact that the mean $\mu$ and autocovariance
function $\chi(\bf r)$ are independent of the variable $\bf x$ is a consequence 
of the translational invariance of the measure. Clearly, 
\begin{equation}
\chi(0)=    \langle Y^2 \rangle -\mu^2
\label{chi-Y0}
\end{equation}
and it follows from Schwarz's inequality that
\begin{equation}
|\chi({\bf r})| \le  \langle Y^2 \rangle -\mu^2.
\label{chi-Y}
\end{equation}
It immediately follows \cite{Lo63} that the autocovariance
function $\chi({\bf r})$ must be positive
semidefinite (nonnegative) in the sense that for any finite number
of spatial locations ${\bf r}_1,{\bf r}_2,\ldots, {\bf r}_m$ in $\mathbb{R}^d$
and arbitrary real numbers $a_1,a_2,\ldots,a_m$,
\begin{equation}
\sum_{i=1}^m \sum_{j=1}^m a_ia_j \chi({\ bf r}_i-{\bf r}_j) \ge 0.
\label{chi-pos}
\end{equation}
It is clear that $\langle [Y({\bf x+ r})-Y({\bf x})]^2 \rangle=2[\chi({\bf 0})-\chi({\bf r})]$.
Thus, if the autocovariance function $\chi({\bf r})$ is continuous at the point $\bf r=0$, the
process $Y(\bf x)$ on $\mathbb{R}^d$ will be {\it mean square continuous}, i.e., 
$ \lim_{\bf r \rightarrow 0} \langle [Y({\bf x+ r})-Y({\bf x})]^2 \rangle=0$ for all $\bf x$. 
Stochastic processes that are continuous in the mean
square sense will be called {\it ordinary}. It is simple to show that
if $\chi({\bf r})$ is continuous at $\bf r=0$, it is continuous for all $\bf r$.

Does every continuous positive semidefinite function $f({\bf r})$  correspond to
a translationally invariant ordinary stochastic process with
a continuous autocovariance $\chi({\bf r})$? The answer is yes, and a proof
is given in the book by Lo{\` e}ve \cite{Lo63}
for stochastic processes in time. Here we state without proof the generalization
to stochastic processes in space.

\begin{theorem}
A continuous function $f(\bf r)$ on $\mathbb{R}^d$ is an autocovariance function
of a translationally invariant ordinary  stochastic process if and only if it is
positive semidefinite.
\end{theorem}

\noindent{Remark:}
\smallskip

\noindent Assuming that $f(\bf r)$ is positive semidefinite, one needs 
to show that there exists a random variable $Y(\bf x)$ on $\mathbb{R}^d$ such 
that $\langle [Y({\bf x})-\mu][Y({\bf x}+{\bf r})-\mu] \rangle=f(\bf r)$.
This is done by demonstrating that there exists a Gaussian (normal) process  
for every autocovariance function \cite{Lo63}. A crucial property
of a Gaussian process is that its full probability distribution
is completely determined by its mean and autocovariance.
\bigskip

The nonnegativity condition (\ref{chi-pos}) is difficult to check. It turns out
that it is easier to establish the existence of an autocovariance function
by appealing to its spectral representation.
We  denote the space of absolutely integrable functions on $\mathbb{R}^d$ by $L^1$. 
The Fourier transform of an $L^1$ function $f : \mathbb{R}^d \rightarrow \Re$
is defined by
\begin{equation}
{\tilde f}({{\bf k}})=\int_{\mathbb{R}^d} f({{\bf r}}) e^{\displaystyle -i{{\bf k}}\cdot{{\bf r}}}\, d{{\bf r}}.
\label{spec}
\end{equation}
If  $f : \mathbb{R}^d \rightarrow \mathbb{R}$ is a radial function, i.e., $f$ depends only
on the modulus $r=|\bf r|$ of the vector $\bf r$, then
its Fourier transform given by
\begin{equation}
{\tilde f}(k) =\left(2\pi\right)^{\frac{d}{2}}\int_{0}^{\infty}r^{d-1}f(r)
\frac{J_{\left(d/2\right)-1}\!\left(kr\right)}{\left(kr\right)^{\left(d/2\right
)-1}} \,dr,
\label{fourier}
\end{equation}
where $k$ is the modulus of the wave vector $\bf k$
and $J_{\nu}(x)$ is the Bessel function of order $\nu$.
The Wiener--Khintchine theorem states that a necessary and sufficient condition for
the existence of  a continuous autocovariance function $\chi({\bf r})$
of a translationally invariant stochastic process $\{Y({\bf x}) : {\bf x}\in \mathbb{R}^d\}$
is that its Fourier transform is nonnegative everywhere, i.e., ${\tilde \chi}({\bf k}) \ge 0$ for all $\bf k$
\cite{Ya87,To02a}. The key ``necessary" part of the proof of this theorem rests on a well-known
theorem due to Bochner \cite{Bo36}, which states that any
continuous function $f({\bf r})$ is positive semidefinite
in the sense of (\ref{chi-pos}) if and only if it has a Fourier--Stieltjes
representation with a nonnegative bounded measure.

\subsection{Generalized Stochastic Processes}

The types of autocovariance functions that we are interested in must allow for
generalized functions, such as Dirac delta functions. The Weiner-Khintchine
theorem has been extended to autocovariances in the class of generalized functions
called {\it tempered distributions}, i.e.,
continuous linear functionals on the space $S$ of infinitely differentiable
functions $\Phi({\bf x})$ such that $\Phi({\bf x})$ as well as all of its
derivatives decay faster than polynomially. Nonnegative
tempered distributions are nonnegative unbounded measures $\nu$ on $\mathbb{R}^d$ such that
\begin{equation}
\int_{\mathbb{R}^d} \frac{d\nu({\bf r})}{(1+|{\bf r}|)^n} < \infty
\end{equation}
for some $n$.
The interested reader is referred to the books by Gel'fand \cite{Ge64} and
Yaglom \cite{Ya87} for details about generalized stochastic
processes. It suffices to say here that 
$\{Y(\Phi({\bf x})) : {\bf x}\in \mathbb{R}^d\}$ is a generalized stochastic process
if for each $\Phi({\bf x}) \in S$ we have a random variable $Y(\Phi({\bf x}))$
that is linear and mean square continuous in $\Phi$.
Then the  mean is the linear functional $\mu(\Phi({\bf x}))=\langle Y(\Phi_1({\bf x})) \rangle$
and the autocovariance function is the bilinear functional 
$\langle [Y(\Phi_1({\bf x}))-\mu(\Phi({\bf x}_1))][Y(\Phi_2({\bf x}+{\bf r}))-
\mu(\Phi({\bf x}_2))] \rangle$, which we still denote by $\chi({\bf r})$ for simplicity.

\begin{theorem}
A necessary and sufficient condition for
an autocovariance function $\chi({\bf r})$
to come from a translationally invariant generalized stochastic process $\{Y(\Phi({\bf x})) : {\bf x}\in \mathbb{R}^d\}$ 
is that that its Fourier transform ${\tilde \chi}({\bf k})$ is a nonnegative 
tempered distribution.
\label{th-wiener}
\end{theorem}
\bigskip

\noindent {\bf Remark:}

\noindent We will call Theorem \ref{th-wiener} the generalized Wiener-Khintchine theorem.

\subsection{Stochastic Point Processes}

Loosely speaking, a stochastic point process  in $\mathbb{R}^d$ 
is defined as a mapping from a probability space
to configurations of points ${\bf x}_1, {\bf x}_2, {\bf x}_3\ldots$
in $\mathbb{R}^d$. More precisely, let $X$ denote
the set of configurations such that each configuration 
$x \in X$ is a subset of $\mathbb{R}^d$ that satisfies two regularity
conditions: (i) there are no multiple points
(${\bf x}_i \neq {\bf x}_j$ if $ i\neq j$) and (ii)
each bounded subset of $\mathbb{R}^d$ must contain
only a finite number of points of $x$.
We denote by $N(B)$  the number of points within
$x \cap B$, $B \in {\cal B}$, where ${\cal B}$ is
the ring of bounded Borel sets in $R^d$. Thus, we always
have $N(B) < \infty$ for $B \in {\cal B}$ but the
possibility $N(\mathbb{R}^d)=\infty$ is not excluded.
We denote by ${\cal U}$ the minimal $\sigma$-algebra of subsets
of $X$ that renders all of the functions $N(B)$ measurable.
Let $(\Omega,{\cal F}, {\cal P})$ be a probability
space. Any measurable map $x(\omega): \Omega \rightarrow X, \omega\in
\Omega$, is called a stochastic point process \cite{St95}.
Point processes belong to the class of generalized stochastic
processes.

A particular realization of a point process 
in $\mathbb{R}^d$ can formally be characterized by the random variable 
\begin{equation}
n({\bf r})=\sum_{i=1}^{\infty} \delta({\bf r} -{\bf x}_i)
\label{n}
\end{equation}
called the ``local" density at position $\bf r$, where $\delta({\bf r})$ 
is a $d$-dimensional Dirac delta function. The ``randomness"
of the point process enters through the positions ${\bf x}_1,{\bf
x}_2,\ldots$. Let us call
\begin{equation}
I_{A}({\bf r})=\Bigg\{{1, \quad {\bf r} \in A,
\atop{0, \quad {\bf r} \notin A,}}
\label{window}
\end{equation}
the {\it indicator function} of a measurable set $A \subset \mathbb{R}^d$, which we 
also call a ``window." For a particular realization,
the number of points $N(A)$ within such a window
is given by
\begin{eqnarray}
N(A) &=& \int_{\mathbb{R}^d} n({\bf r})
I_{A}({\bf r}) \, d{\bf r} \nonumber \\
&=& \sum_{i=1}^{\infty} \int_{\mathbb{R}^d} \delta({\bf r}-{\bf x}_i) I_A({\bf r}) \, d{\bf r}
\nonumber \\
&=& \sum_{i\ge 1} I_A({\bf x}_i).
\label{N}
\end{eqnarray}
Note that this random setting is quite general.
It incorporates cases in which the location of the points are deterministically
known, such as a lattice. A packing of congruent spheres of unit diameter
is simply a point process in which any pair of points cannot be closer
than a unit distance from one another.

It is known that the probability measure on $(X,{\cal U})$ exists provided that
the infinite set of $n$-point correlation functions $\rho_n$,
$n=1,2,3\ldots$ meet certain conditions \cite{Le73,Le75a,Le75b}.
The $n$-point correlation function $\rho_n({\bf r}_1,{\bf r}_2,\ldots,{\bf r}_n)$
is the contribution to the expectation $\langle n({\bf r}_1) n({\bf r}_2)
\cdots n({\bf r}_n)\rangle$ when the indices on the sums are not equal to one another, i.e.,
\begin{equation}
\rho_n({\bf r}_1,{\bf r}_2,\ldots,{\bf r}_n)= \left\langle \sum_{i_1\neq i_2 \neq \cdots \neq i_n}^{\infty}  
\delta({\bf r}_1 -{\bf x}_{i_1}) \delta({\bf r}_2 -{\bf x}_{i_2})\cdots \delta({\bf r}_n -{\bf x}_{i_n})\right\rangle.
\end{equation}
Note that the distribution-valued function $\rho_n({\bf r}_1,{\bf r}_2,\ldots,{\bf r}_n)$ also
has a probabilistic interpretation: apart from trivial constants,
it is the probability density function
associated with finding $n$ different points at positions 
${\bf r}_1,{\bf r_2},\ldots,{\bf r}_n$. For this reason, $\rho_n$ 
is also called the $n$-particle density and, for any $n$,
has the nonnegativity property
\begin{equation}
\rho_n({\bf r}_1,{\bf r}_2,\ldots,{\bf r}_n) \ge 0 \qquad \forall {\bf r}_i \in \mathbb{R}^d \quad (i=1,2,\ldots n).
\label{positive}
\end{equation}
Translational invariance means that for every constant vector $\bf y$ in $\mathbb{R}^d$, 
$\rho_n({\bf r}_1,{\bf r}_2,\ldots,{\bf r}_n)=\rho_n({\bf r}_1+{\bf y},\ldots,{\bf r}_n+{\bf y})$,
which implies that   
\begin{equation}
\rho_n({\bf r}_1,{\bf r}_2,\ldots,{\bf r}_n)=\rho^ng_n({\bf r}_{12},\ldots, {\bf r}_{1n}),
\label{nbody}
\end{equation}
where $\rho$ is the number (or center) density
and $g_n({\bf r}_{12},\ldots, {\bf r}_{1n})$ is the {\it $n$-particle correlation function},
which depends on the relative positions ${\bf r}_{12}, {\bf r}_{13}, \ldots$,
where ${\bf r}_{ij} \equiv {\bf r}_j -{\bf r}_i$ and we have chosen the origin to be at ${\bf r}_1$.

For such point processes without {\it long-range order}, 
$g_n({\bf r}_{12},\ldots, {\bf r}_{1n}) \rightarrow 1$ when
the points (or ``particles") are mutually far from one another, i.e.,  
as $|{\bf r}_{ij}| \rightarrow\infty$ 
($1\leq i < j < \infty$), $\rho_n({\bf r}_{1}, {\bf r}_2, \ldots, {\bf r}_{n}) \rightarrow \rho^n$.
Thus, the deviation of $g_n$ from unity  provides a
measure of the degree of spatial correlation
between the particles, with unity corresponding to no spatial correlation.
Note that for a translationally invariant Poisson point process, 
$g_n$ is unity for all values of its argument. 

As we indicated in the beginning of this section,
the first two correlation functions, $\rho_1({\bf r}_1)=\rho$ and $\rho_2({\bf r}_1,{\bf r}_2)=\rho^2 g_2({\bf r})$, 
for translationally invariant point processes are of central concern in this paper.
%When such a point process does not possess {\it long-range order}, i.e., $\rho_2({\bf r}) \rightarrow \rho^2$
%as $|\bf r| \rightarrow \infty$, it is convenient to define the {\it pair correlation function} 
%\begin{equation}
%g_2({\bf r}) \equiv \frac{\rho_2({\bf r})}{\rho^2}.
%\label{g2-rho2}
%\end{equation}
%Thus,  $g_2({\bf r})$ tends to unity for large $|\bf r|$ in such instances,  indicating that
%there is no spatial correlation
%between a pair of widely separated points. (Note that for a Poisson
%point process, $g_2$ is unity for all pair distances.)
If the point process is also rotationally invariant
(statistically isotropic), then 
$g_2$ depends on  the radial distance $r=|\bf r|$ only, i.e.,
\begin{equation}
g_2({\bf r}) = g_2(r),
\end{equation}
and is referred to as the {\it radial distribution function}. Strictly speaking,
one should use different notation for the left and right members
of Eq. (22), but to conform to conventional statistical-mechanical usage,
we invoke the common notation for both.
Because $\rho_2({\bf r}_1,{\bf r}_2)/\rho= \rho g_2(r)$ is a conditional
joint probability density, then 
\begin{equation}
Z(r_1,r_2)= \int_{r_1}^{r_2} \rho s_1(r) g_2(r) \, dr
\end{equation}
is the expected number of points at radial distances between $r_1$ and $r_2$
from a randomly chosen point.
Here $s_1(r)$ is the surface area of a $d$-dimensional sphere of radius $r$ given by
\begin{equation}
s_1(r)=\frac{2\pi^{d/2}r^{d-1}}{\Gamma(d/2)}.
\label{surf}
\end{equation}
For a packing of congruent spheres of unit diameter, $g(r)=0$ for $0 \le r <1$, i.e.,
\begin{equation}
\mbox{supp}(g_2) \subseteq \{r: r \ge 1\}.
\label{support}
\end{equation}
Note that the radial distribution function $g_2(r)$ (or any of the $\rho_n$)
for a point process must be able to incorporate Dirac delta functions. We will
specifically consider those radial distribution functions that are
nonnegative distributions. For example, $g_2(r)$
for a lattice packing is the rotational symmetrization of
the sum of delta functions at lattice points at a radial
distance $r$ from any lattice point \cite{To03a}.

For a translationally invariant point process, the autocovariance function
$\chi({\bf r})$ is related to the pair correlation function via
\begin{equation}
\chi({\bf r})=\rho \delta({\bf r})+ \rho^2 h({\bf r}),
\end{equation}
where
\begin{equation}
h({\bf r})\equiv g_2({\bf r})-1,
\label{total}
\end{equation}
is the {\it total correlation function.} This relation is obtained using definitions
(\ref{auto}) and (\ref{n}) with $Y({\bf x})=n({\bf x})$. Note that $\chi({\bf r})=\rho
\delta({\bf r})$ (i.e., $h=0$) for a translationally invariant Poisson point process.
{\it Positive and negative pair correlations} are defined as those instances in which $h$ is
positive (i.e., $g_2 >1$) and $h$ is negative (i.e., $g_2<1$), respectively.
The Fourier transform of the distribution-valued function $\chi({\bf r})$ is given by
\begin{equation}
{\tilde \chi}({\bf k})=\rho+\rho^2 {\tilde h}({\bf k}),
\label{chi2}
\end{equation}
where ${\tilde h}({\bf k})$ is the Fourier transform of $h(\bf r)$. It is common
practice in statistical physics to deal with a function trivially related
to the spectral function ${\tilde \chi}({\bf k})$ called the structure factor 
$S({\bf k})$, i.e.,
\begin{equation}
S({\bf k}) \equiv \frac{{\tilde \chi}({\bf k})}{\rho}=1+\rho{\tilde h}({\bf k}).
\label{factor}
\end{equation}

A natural question to ask at this point is the following: Given a positive
number density $\rho$ and a  pair correlation function $g_2({\bf r})$, does
there exist a translationally invariant point process in $\mathbb{R}^d$ with measure $\cal P$
for which $\rho$ and $g_2$ are one-point and two-point correlation functions,
respectively. Two obvious nonnegativity
conditions \cite{To02c} that must be satisfied are the following:
\begin{equation}
g_2({\bf r}) \ge 0 \qquad \mbox{for all}\quad {\bf r}
\label{cond1}
\end{equation}
and 
\begin{equation}
S({\bf k})=1+\rho {\tilde h}({\bf k}) \ge 0  \qquad \mbox{for all}\quad {\bf k}.
\label{cond2}
\end{equation}
The first condition is trivial and comes  from (\ref{positive}) with $n=2$. The second condition
is nontrivial and derives from the  generalized Wiener-Khintchine
theorem (\ref{th-wiener}) using relations ({\ref{chi2}) and ({\ref{factor}).
However, for realizability of point processes in arbitrary dimension $d$, the
two standard conditions (\ref{cond1}) and (\ref{cond2})
are only necessary, not necessary and sufficient. The same state of affairs applies to the theory of random sets
\cite{To02a}, where it is known that the Wiener-Khintchine theorem only
provides a necessary condition on realizable autocovariance functions.
The simplest example of a random set is one in which $\mathbb{R}^d$
is partitioned into two disjoint regions (phases) but with an interface
that is only known in a probabilistic sense. (A packing can therefore
be viewed as a special random set). Thus, a random set
is described by a random variable that is the
indicator function for a particular phase, i.e., it
is a binary stochastic process. The class  of autocovariances
that comes from a binary stochastic process
is a subclass of the total class that comes
from an ordinary process $\{Y({\bf x}) : {\bf x}\in \mathbb{R}^d\}$
and meets the existence condition of Theorem~\ref{th-wiener}. 
Similarly, the class of autocovariances
that comes from a  point process is a subset of 
of a generalized process $\{Y(\Phi({\bf x})) : {\bf x}\in \mathbb{R}^d\}$
and therefore the  existence condition of Theorem~\ref{th-wiener} is
only necessary. 

It has recently come to light that a positive $g_2$ for a positive $\rho$ must
satisfy an uncountable number of necessary and sufficient
conditions for it to correspond to a  realizable  point process \cite{Cos04}.
However, these conditions are very difficult (or, more likely, impossible)
to check for arbitrary dimension. In other words, given $\rho_1=\rho$
and $\rho_2=\rho^2 g_2$, it is difficult to ascertain if there
are some higher-order functions $\rho_3, \rho_4, \ldots$ for which these one- and two-point
correlation functions hold. Thus, an important practical problem
becomes the determination of a manageable number of
necessary conditions that can be readily checked. 

One such additional necessary condition, obtained by Yamada \cite{Ya61},
is concerned with $\sigma^2(A) \equiv 
\langle (N(A)- \langle N(A) \rangle)^2\rangle$,
the variance in the number of points $N(A)$ contained within a window
$A \subset \mathbb{R}^d$. Specifically, he showed that
\begin{equation}
\sigma^2(A) \ge \theta(1-\theta),
\label{yamada}
\end{equation}
where $\theta$ is the fractional part of the expected number of points 
$\rho |A|$ contained in
the window. This inequality is a consequence of the fact that the number
of points $N(A)$ within a window at some fixed position is an integer, not a continuous
variable, and sets a lower limit on the number variance. We note in passing
that the determination of the number variance for lattice point patterns is an outstanding
problem in number theory \cite{Ke48,Ke53b,Sa05}. The number variance for a specific choice
of $A$ is necessarily a positive number and 
generally related to the total pair correlation function $h(\bf r)$ for
a translationally invariant point process \cite{To03a}.
In the special case of a spherical window of radius $R$ in $\mathbb{R}^d$,
it is explicitly given by
\begin{equation}
\sigma^2(R)=\rho v_1(R) \Bigg[ 1+\rho \int_{\mathbb{R}^d} h({\bf r}) \alpha_2({\bf r}; R) \, d{\bf r}\Bigg] \ge \theta(1-\theta),
\label{variance}
\end{equation}
where $\sigma^2(R)$ is the number variance for
a spherical window of radius $R$ in $\mathbb{R}^d$, $v_1(R)$ is the volume of the window,
and $\alpha_2(r;R)$ is the volume common to two spherical windows of radius $R$
whose centers are separated by a distance $r$ divided by $v_1(R)$.
We will call $\alpha_2(r;R)$ the {\it scaled intersection volume}.
The lower bound (\ref{variance}) provides another integral 
condition on  the pair correlation function.

For large $R$, it has been proved that $\sigma^2(R)$ cannot grow more
slowly than $\gamma R^{d-1}$, where $\gamma$ is a positive constant \cite{Beck87}.
This implies that the Yamada lower bound in (\ref{variance}) is always satisfied for 
for sufficiently large $R$ for any $d \ge 2$. In fact, we have not been
able to construct any examples of a pair correlation function $g_2({\bf r})$
at some number density $\rho$ that satisfy the two nonnegativity
conditions (\ref{cond1}) and (\ref{cond2}) and  simultaneously
violate the Yamada condition for any $R$ and any $d \ge 2$.
Thus, it appears that the Yamada condition is most relevant
in one dimension, especially in those cases when $\sigma^2(R)$
is bounded. We note that point processes (translationally invariant or not)
for which $\sigma^2(R) \sim R^{d-1}$ for large $R$ are examples
of {\it hyperuniform} point patterns \cite{To03a}. This classification
includes all lattices as well as aperiodic point patterns. Hyperuniformity
implies that the structure factor $S({\bf k})$ has the following
small ${\bf k}$ behavior:
\begin{equation}
\lim_{{\bf k}\rightarrow {\bf 0}} S({\bf k}) =0.
\label{hyper}
\end{equation}

The scaled intersection volume  $\alpha_2(r;R)$ appearing in (\ref{variance}) will play a prominent role
in this paper. It has  support in the
interval $[0,2R)$, the range $[0,1]$ and the following  integral representation:
\begin{equation}
\alpha_2(r;R) = c(d) \int_0^{\cos^{-1}(r/(2R))} \sin^d(\theta) \, d\theta,
\label{alpha}
\end{equation}
where $c(d)$ is the $d$-dimensional constant given by
\begin{equation}
c(d)= \frac{2 \Gamma(1+d/2)}{\pi^{1/2} \Gamma((d+1)/2)}.
\end{equation}
Following the analysis given by Torquato and Stillinger \cite{To03a} for low
dimensions, we obtain the following {\it new} series representation
of the scaled intersection volume $\alpha_2(r;R)$ for $r \le 2R$
and for any $d$:
\begin{equation}
\alpha(r;R)=1- c(d) x+c(d) \sum_{n=2}^{\infty}
(-1)^n \frac{(d-1)(d-3) \cdots (d-2n+5)(d-2n+3)}
{(2n-1) 2 \cdot 4 \cdot 6 \cdot (2n-4)(2n-2)} x^{2n-1},
\label{series}
\end{equation}
where $x=r/(2R)$. This is easily proved with the help of Maple. For even dimensions, relation (\ref{series}) is
an infinite series, but for odd dimensions, the series truncates such
that $\alpha_2(r;R)$ is a univariate polynomial of degree $d$. Except for
the first term of unity, all of the terms in (\ref{series}) involve only odd
powers of $x$. Figure \ref{intersection} shows graphs of the scaled intersection
volume $\alpha_2(r;R)$ as a function of $r$ for the first
five space dimensions. For any dimension, $\alpha(r;R)$
is a monotonically decreasing function of $r$. At a fixed
value of $r$ in the interval $(0,2R)$, $\alpha_2(r;R)$
is a monotonically decreasing function of the dimension $d$.
For large $d$, we will subsequently make use of the asymptotic result
\begin{equation}
\alpha_2(R;R) \sim \left(\frac{6}{\pi}\right)^{1/2}  \left(\frac{3}{4}\right)^{d/2}
\frac{1}{d^{1/2}}.
\label{alpha-asym}
\end{equation}

\begin{figure}[bthp]
\centerline{\psfig{file=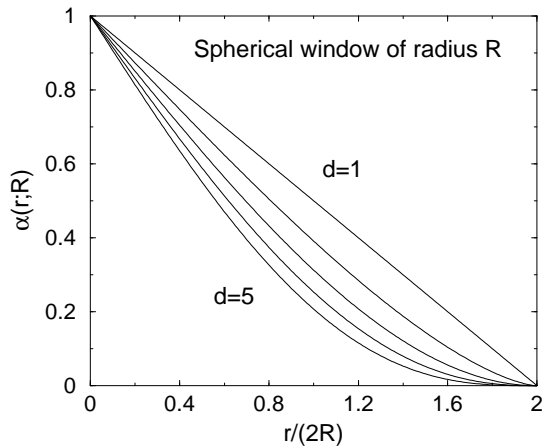,width=2.8in}}
\caption{The scaled intersection
volume $\alpha_2(r;R)$ for spherical windows of radius $R$
as a function of $r$ for the first
five space dimensions. The uppermost curve is for $d=1$
and lowermost curve is for $d=5$.}
\label{intersection}
\end{figure}

Before closing this section, it is useful
to note that there has been some recent work
that demonstrates the existence of point processes for 
a specific  $\rho$ and $g_2$ provided that
$\rho$ and $g_2$ meet certain restrictions. 
For example, Ambartzumian and Sukiasian proved the existence
of point processes that come from Gibbs measures for a special $g_2$ for sufficiently
small $\rho$ \cite{Am91}. Determinantal point processes have been
considered by Soshnikov \cite{So00} and Costin and Lebowitz \cite{Cos04}. Costin and Lebowitz
have also studied certain one-dimensional renewal point processes \cite{Cos04}.
Stillinger and Torquato \cite{St04} discussed the possible existence
of a general interparticle pair potential (associated with a Gibbs measure) for a given
$\rho$ and $g_2$ using a cluster expansion procedure but did not address the
issue of convergence of this expansion. Koralov \cite{Ko05} indeed
proves the existence of a pair potential on a lattice (with
the restriction of single occupancy per lattice site) for
which $\rho$ is the density and $g_2$ is the pair correlation
function for sufficiently small $\rho$ and $g_2$. There is no
reason to believe that Koralov's proof is not directly extendable to the 
case of a point process corresponding to a sphere
packing in $\mathbb{R}^d$, where the nonoverlap condition is
the analog of single occupancy on the lattice.
Thus, we expect that one can prove the existence of a pair potential in $\mathbb{R}^d$
corresponding to a sphere packing for a given $\rho$ and $g_2$ 
provided that $\rho$ and $g_2$ sufficiently small.

\section{Disordered Packings in High Dimensions and the Decorrelation Principle}
\label{disordered}

In this section,  we examine the asymptotic behavior
of certain disordered packings in high dimensions and show
that unconstrained spatial correlations vanish asymptotically, yielding
a {\it decorrelation principle}. We define a {\it disordered packing} in $\mathbb{R}^d$
to be one in which the pair correlation function $g_2({\bf r})$ decays
to its long-range value of unity faster than $|{\bf r}|^{-d-\varepsilon}$ for
some $\varepsilon >0$.
The decorrelation principle as well as a number of other results (that will be
discussed in Section 5)  motivate us to propose a conjecture
in Section 5 that describes the circumstances in which the two standard
nonnegativity conditions (\ref{cond1}) and (\ref{cond2})
are necessary and sufficient to ensure the existence
of a disordered sphere packing.

\subsection{Example 1: Disordered Sequential Packings}

First we show that there exists a disordered sphere packing that realizes the greedy lower bound
(\ref{sat}) ($\phi = 1/2^d$) for all $d$. 
Then we study the asymptotic properties of the $n$-particle correlation functions
in the large dimension limit.

The disordered packing that achieves the greedy lower bound is a special
case of a generalization of the so-called random sequential addition (RSA) process \cite{To02a}.
This generalization, which we introduce here, is a subset of the Poisson point process. Specifically, 
the centers of ``test" spheres of unit diameter arrive continually throughout $\mathbb{R}^d$
during time $t\ge 0$ according to a translationally invariant Poisson process
of density per unit time $\eta$, i.e., $\eta$ is the number of
points per unit volume and time. Therefore, the expected number of
centers in a region of volume $\Omega$ during time $t$ is $\eta \Omega t$
and the probability that this region is empty of centers is $\exp(-\eta \Omega t)$.
However, this Poisson distribution of test spheres is not a packing because 
the spheres can overlap. To create a packing from this point process,
one must remove test spheres such that no sphere center can lie within a spherical 
region of unit radius from any sphere center. Without loss of generality,
we will set $\eta=1$.

There is a variety of ways of achieving this ``thinning" process such that the subset of points
correspond to a sphere packing. One obvious rule is to retain a test
sphere at time $t$ only if it does not overlap a sphere that was successfully
added to the packing at an earlier time. This criterion defines the well-known RSA process in $\mathbb{R}^d$
\cite{To02a}, and is clearly a statistically homogeneous and isotropic
sphere packing in $\mathbb{R}^d$ with a time-dependent density
$\phi(t)$. In the limit $t \rightarrow \infty$, the RSA process corresponds to a saturated
packing with a maximal or {\it saturation} density $\phi_s(\infty) \equiv \lim_{t\rightarrow
\infty} \phi(t)$. In one dimension, the RSA process is commonly known as the ``car parking problem", 
which Re{\' n}yi showed has a saturation density $\phi_s(\infty)= 0.7476\ldots$ \cite{Re63}.
For $2 \le d < \infty$, an exact determination of $\phi_s(\infty)$ is not
possible, but estimates for it have been obtained via computer experiments
for low dimensions \cite{To02a}. However, as we will discuss below, the standard
RSA process at small times (or, equivalently, small densities) can be analyzed exactly.

Another thinning criterion retains a test sphere centered at position $\bf r$ at time $t$ 
if no other test sphere is within a unit radial distance from $\bf r$ for
the time interval $\kappa t$ prior to $t$, where $\kappa$ is a positive constant
in the interval $[0,1]$. This packing is a subset of the RSA packing, and
therefore we refer to it as the generalized RSA process. Note
that when $\kappa=0$, the standard RSA process is recovered, and when $\kappa=1$,
a relatively unknown model due to Mat{\' e}rn \cite{Ma86} is recovered. The latter is amenable
to exact analysis. 

The time-dependent density $\phi(t)$ in the case of the generalized RSA process
with $\kappa=1$ is easily obtained. (Note that for any $0 < \kappa \le 1$, the generalized
RSA process is always an {\it unsaturated} packing.) In this packing,
a test sphere at time $t$ is accepted only if does not overlap 
an existing sphere in the packing as well as any previously rejected
test sphere (which we will call ``ghost" spheres.) An overlap cannot
occur if a test sphere is outside a unit radius of any successfully
added sphere or ghost sphere. Because of the underlying Poisson process,
the probability that a trial sphere is retained
at time $t$ is given by $\exp(-v_1(1) t)$, where $v_1(1)$ is the volume
of a sphere of unit radius having the same center
as the retained sphere of radius 1/2. Therefore, the expected number density
$\rho(t)$ and packing density $\phi(t)$ at any time $t$ are respectively  given by
\begin{equation}
\rho(t)=\int_0^t \exp(-v_1(1) t^\prime) \, dt^\prime=\frac{1}{v_1(1)}[1-\exp(-v_1(1) t)]
\end{equation}
and
\begin{equation}
\phi(t)= \rho(t) v_1(1/2)=\frac{1}{2^d}[1-\exp(-v_1(1) t)].
\label{phi}
\end{equation}
We see that $\phi(t)$ is a monotonically increasing function of $t$.
This result was first given by Mat{\' e}rn using a different approach
and he also gave a formal expression for
the time-dependent radial distribution function $g_2(r;t)$ (see Section \ref{realize}). 
Here we present an explicit expression for $g_2(r;t)$ at time $t$ for any dimension $d$:
\begin{equation}
g_2(r;t)=\frac{\Theta(r-1)}{2^{2d-1}[\beta_2(r;1)-1]\phi^2(t)}
\Bigg[ 2^d \phi(t)- \frac{1-\exp[-2^d \beta_2(r;1)t]}{\beta_2(r;1)}\Bigg].
\label{g}
\end{equation}
Here 
\begin{equation}
\Theta(x) =\Bigg\{{0, \quad x<0,\atop{1, \quad x \ge 0,}}
\label{heaviside}
\end{equation}
is the unit step function and 
\begin{equation}
\beta_2(r;R)=2- \alpha_2(r;R)
\end{equation}
is the union volume of two spheres of radius $R$ (whose centers are separated by
a distance $r$) divided by the volume of a sphere of radius $R$ and $\alpha_2(r;R)$
is the scaled intersection volume of two such spheres given by (\ref{alpha}).
Our approach for obtaining (\ref{g}) is different than Mat{\' e}rn's
and details are given elsewhere \cite{To05a}.

It is useful to note that at small times or, equivalently, low densities, 
formula (\ref{phi}) yields the asymptotic expansion
$\phi(t)=v_1(1)t/2^d-  v^2_1(1)t^2/2^{d+1}+{\cal O}(t^3)$,
which when inverted yields $t=2^d\phi/v_1(1)+ 2^{d-1} \phi^2+{\cal O}(\phi^3)$.
Substitution of this last result into (\ref{g}) gives
\begin{equation}
g_2(r;\phi)=\Theta(r-1)+ {\cal O}(\phi^3),
\label{g2-t-grsa}
\end{equation}
which implies that $g_2(r;\phi)$ tends to the unit step function
$\Theta(r-1)$ as $\phi \rightarrow 0$ for any finite $d$.
\begin{figure}[bthp]
\centerline{\psfig{file=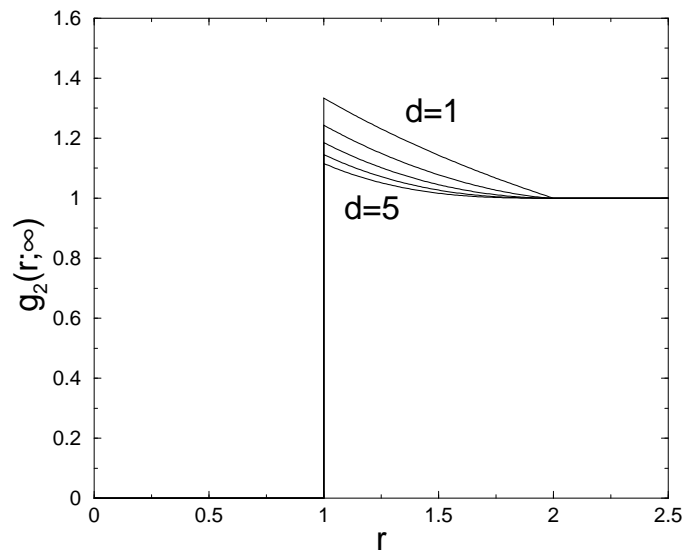,height=3.0in}}
\caption{ Radial distribution function for the first five space dimensions
at the maximum density $\phi=1/2^d$ for the generalized RSA model with $\kappa=1$.}
\label{grsa}
\end{figure}

In the limit $t \rightarrow \infty$, the maximum density
is given by
\begin{equation}
\phi(\infty) \equiv \lim_{t \rightarrow \infty} \phi(t)=\frac{1}{2^d}
\end{equation}
and 
\begin{equation}
g_2(r;\infty) \equiv \lim_{t \rightarrow \infty} g_2(r;t)=\frac{2\Theta(r-1)}{\beta_2(r;1)}=
\frac{\Theta(r-1)}{1-\alpha_2(r;1)/2}.
\label{g2-grsa}
\end{equation}
We see that the greedy lower-bound limit on the density is achieved in 
the infinite-time limit for this sequential but unsaturated packing. This
is the first time that such an observation has been made. Obviously, for
any $0 \le \kappa < 1$, the maximum (infinite-time) density 
of the generalized RSA packing is bounded from below by $1/2^d$ (the
maximum density for $\kappa=1$).
Note also that because $\beta_2(r;1)$ is equal to 2 for $r \ge 2$,  
$g_2(r;\infty)=1$ for $r \ge 2$, i.e., spatial correlations
vanish identically for all pair distances except
those in the small interval $[0,2)$. Even the positive correlations
exhibited for $1 < r <2$ are rather weak and decrease 
with increasing dimension. The function $g_2(r;\infty)$ achieves its largest value
at $r=1^+$ in any dimension and for $d=1$, $g_2(1^+;\infty)=4/3$.
The radial distribution function $g_2(r;\infty)$ is plotted in Fig. \ref{grsa}
for the first five space dimensions. Using the asymptotic result (\ref{alpha-asym})
and relation (\ref{g2-grsa}),
it is seen that for large $d$,
\begin{equation}
g_2(1^{+};\infty) \sim \frac{ \Theta(r-1)}{\displaystyle 1- \left(\frac{3}{2\pi}\right)^{1/2}  \left(\frac{3}{4}\right)^{d/2}
\frac{1}{d^{1/2}}},
\end{equation}
and thus $g_2(r;\infty)$ tends  to the unit step function
$\Theta(r-1)$ exponentially fast as $d \rightarrow \infty$ because the scaled
intersection volume $\alpha_2(1;1)$ vanishes exponentially fast.

The higher-order correlation functions for this model have not
been given previously. In another work \cite{To05a}, we use an approach different
from the one used by Mat{\' e}rn to obtain not only $g_2$ but an explicit formula
 for the general $n$-particle correlation function $g_n$, defined by (\ref{nbody}), for any time
$t$ and $n$ and for arbitrary dimension $d$. To our
knowledge, this represents the first
exactly solvable disordered sphere packing
model for any $d$. These details are somewhat tangential to the present work
and for our purposes it suffices to state the final result
in the limit $t \rightarrow \infty$ for $n \ge 2$:
\begin{equation}
g_n({\bf r}_{1}, \ldots, {\bf r}_{n};\infty)= 
\frac{\displaystyle \prod_{i<j}^n\Theta(r_{ij}-1)}{
\beta_n({\bf r}_{1},\ldots,{\bf r}_{n};1)}\Big[\sum_{i=1}^{n}g_{n-1}(Q_i;\infty)\Big]
\label{gn-grsa}
\end{equation}
where the sum is over  all the $n$ distinguishable ways of choosing
$n-1$ positions from $n$ positions ${\bf r}_1,\ldots {\bf r}_n$ and the arguments of $g_{n-1}$
are the associated $n-1$ positions, which we denote by $Q_i$, and $g_1\equiv 1$.
Moreover, $\beta_n({\bf r}_{1},\ldots, {\bf r}_{n};R)$ is the union volume
of $n$ congruent spheres of radius $R$, whose centers are located at
${\bf r}_1,\ldots, {\bf r_n}$, where $r_{ij}=|{\bf r}_j -{\bf r}_i|$
for all $1 \le i <j \le n$, divided by the volume of a sphere of radius $R$.

\begin{lemma}
In the limit $d \rightarrow \infty$, the $n$-particle correlation function 
$g_n({\bf r}_{1}, \ldots, {\bf r}_{n};\infty) \sim 1$ 
uniformly in $({\bf r}_{1}, \ldots, {\bf r}_{n})$ $\in \mathbb{R}^d$
such that $r_{ij}\ge 1$ for all $1 \le i <j \le n$. If  $r_{ij} < 1$ for any 
pair of points ${\bf r}_i$ and ${\bf r}_j$, then $g_n({\bf r}_{1}, \ldots, {\bf r}_{n};\infty)=0$.
\end{lemma}

\noindent
{\it Proof}: The second part of the Lemma is the trivial requirement for a packing.
Whenever $r_{ij}\ge 1$ for all $1 \le i <j \le n$, it is clear from 
(\ref{gn-grsa}) that we have the following upper and lower bounds
on the $n$-particle correlation function:
\begin{equation}
\frac{n}{\beta_n} \le g_n \le \frac{n g_{n-1}^*}{\beta_n},
\end{equation}
where $g_{n-1}^*$ denotes the largest possible value of $g_{n-1}$. The scaled
union volume $\beta_n$ of $n$ spheres obeys the bounds
\begin{equation}
n -\sum_{i<j} \alpha_2(r_{ij};1) \le \beta_n \le n,
\end{equation}
but since the scaled intersection volume of two spheres $\alpha_2(r;1)$
attains its maximum value for $r \ge 1$ when $r=1$, we also have
\begin{equation}
n - \frac{n(n-1)}{2}\alpha_2(1;1) \le \beta_n \le n.
\end{equation}
Use of this inequality and the recursive relation (\ref{gn-grsa}) yields
the bounds
\begin{equation}
1\le g_n \le \frac{1}{\displaystyle 1-\frac{n(n-1) \alpha_2(1;1)}{4}+ {\cal O}(\alpha_2(1;1)^2)}.
\end{equation}
Using the asymptotic result (\ref{alpha-asym}), we see that the upper bound
tends to the lower bound for any given $n$ as $d\rightarrow \infty$, which proves the Lemma.

In summary, the Lemma enables us to conclude that in the limit $d \rightarrow \infty$ 
and for $\phi=1/2^d$
\begin{equation}
g_n({\bf r}_{12}, \ldots, {\bf r}_{1n};\infty) \sim \prod_{i<j}^n g_2(r_{ij};\infty),
\end{equation}
where
\begin{equation}
g_2(r;\infty) \sim \Theta(r-1).
\label{g2-asym-grsa}
\end{equation}

Importantly, we see that the asymptotic behavior
of $g_2$ in the low-density limit $\phi \rightarrow 0$ for any $d$ [cf. (\ref{g2-t-grsa})] is the same as
the high-dimensional limit $d \rightarrow \infty$ [cf. (\ref{g2-asym-grsa})], i.e., spatial
correlations, which exist for positive densities at fixed $d$, vanish
for pair distances beyond the hard-core diameter. Note also that
$g_n$ for $n \ge 3$ asymptotically factorizes into products involving only the
pair correlation function $g_2$. Is the similarity
between the low-density and high-dimensional 
limits for this model of a disordered packing a general characteristic of disordered packings?
In what follows, we discuss another disordered packing that
has this attribute and subsequently formulate what we refer to as
a ``decorrelation principle."

\subsection{Example 2: The Classic Gibbsian Hard-Sphere Packing}

The  statistical mechanics of the classic Gibbsian hard-sphere packing is well
established (see \cite{To02a} and references therein). The purpose of this subsection is simply to
collect some results that motivate the decorrelation principle.
Let $\Phi_N({\bf r}^N)$ be the $N$-body interaction potential 
for a finite but large number of particles with configuration ${\bf r}^N\equiv
\{{\bf r}_{1}, {\bf r}_{2},\ldots, {\bf r}_N\}$ in a volume $V$ in $\mathbb{R}^d$
at absolute temperature $T$. A large collection of such systems in which
$N$, $V$ and $T$ are fixed but in which the particle configurations
are otherwise free to  vary is called the {\it Gibbs canonical ensemble}.
Our interest is in the {\it thermodynamic limit}, i.e., the distinguished limit
in which $N \rightarrow \infty$ and $V \rightarrow \infty$ such
that the number density $\rho=N/V$ exists. For a Gibbs canonical ensemble,
when the $n$-particle densities $\rho_n$ (defined in Section \ref{realize}) exist,
they are entirely determined by the interaction potential
$\Phi_N({\bf r}^N)$. For a hard-sphere packing, the
interaction potential is given by a sum of pairwise terms such that
\begin{equation}
\Phi_N({\bf r}^N) =  \sum_{i<j}^N u_2 (|{\bf r}_j -{\bf r}_i|).
\label{pairwise}
\end{equation}
where $u_2(r)$ is the pair potential defined by
\begin{equation}
u_2(r) =
\left\{\begin{array}{ll}
+ ~\infty, \qquad & r < 1, \nonumber \\
0, \qquad & r \ge 1,\\
\end{array}
\right.
\label{hs-potential}
\end{equation}
Thus,  the particles do not interact for interparticle separation distances
greater than or equal to unity but experience an infinite repulsive force for 
distances less than unity. The hard spheres have kinetic energy, and therefore
a temperature, but the temperature enters in a trivial way
because the configurational correlations between the spheres
are independent of the temperature. We call this the {\it classic  equilibrium}
sphere packing, which is both translationally and rotationally invariant.

In one dimension, the $n$-particle densities $\rho_n$ for such packings are known
exactly in the thermodynamic limit. The density $\phi$ lies in the interval
$[0,1]$ but this one-dimensional packing is devoid of a discontinuous (first-order) transition
from a disordered (liquid) phase to an ordered (solid) phase.
Although a rigorous proof for the existence of a 
liquid-to-solid phase transition in two or three dimensions is not yet available,
there is overwhelming numerical evidence (as obtained from computer simulations) 
that such a transformation takes place at sufficiently high densities.
The maximal densities for equilibrium sphere packings 
in two and three dimensions are $\phi_{\mbox{\scriptsize max}}= \pi/\sqrt{12}$ and $\phi_{\mbox{\scriptsize max}}=\pi/\sqrt{18}$,
respectively, i.e., they correspond to the density of the densest sphere 
packing in the respective dimension. 

\begin{figure}
\centerline{\psfig{file=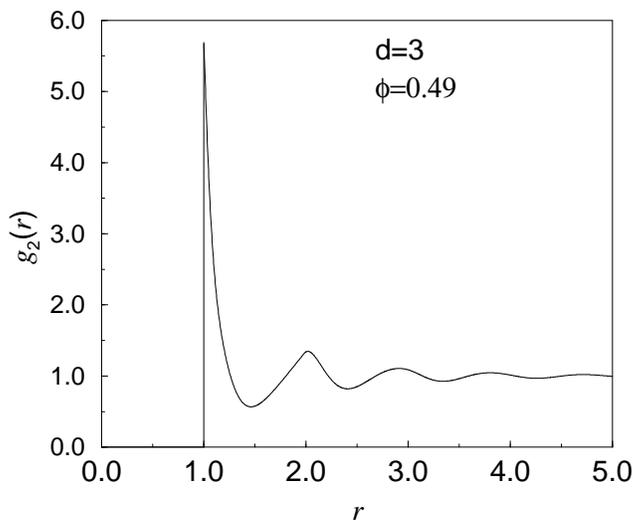,height=2.8in}}
\caption{The radial distribution function
for the classic three-dimensional  equilibrium packing at $\phi=0.49$
as obtained from molecular dynamics computer simulations.
The graph is adapted from Figure 3.15 of Torquato \protect\cite{To02a}.}
\label{rdf}
\end{figure}

Figure \ref{rdf} shows the three-dimensional radial distribution function as obtained from
computer simulations for a density $\phi=0.49$, which is near the maximum value for the 
stable disordered branch. It is seen that the packing exhibits {\it short-range} order
[i.e., $g_2(r)$ has both positive and negative correlations for small $r$],
but $g_2(r)$ decays to its long-range value exponentially fast
after several diameters. By contrast, in the limit $d \rightarrow \infty$,
it has been  shown that the ``pressure" \cite{Ru99} of an equilibrium packing
is exactly given by the first two terms of its asymptotic
low-density expansion for some positive density interval $[0,\phi_0]$ \cite{Wy87,Fr99}.
(Roughly speaking, the pressure is the average force per unit area acting on an ``imaginary
planar wall" in the packing due collisions between the spheres and the wall.)
Frisch and Percus \cite{Fr99} have established, albeit not rigorously,  that
$\phi_0=1/2^d$. This result for the pressure
implies that the leading-order term 
of the low-density expansion of the radial distribution function in arbitrary dimension \cite{To02a}
\begin{equation}
g_2(r)=\Theta(r-1)\Big[1+ 2^d \alpha_2(r;1) \phi + {\cal O}(\phi^2)\Big]
\label{g2-exp}
\end{equation}
becomes asymptotically exact in the limit $d \rightarrow \infty$ in the same density interval.
The presence of the unit step function $\Theta(r-1)$ in relation (\ref{g2-exp})
means that the scaled intersection volume $\alpha_2(r;1)$ need only be considered
for values of $r$ in the interval $[1,2]$. Since $\alpha_2(r;1)$ is largest
when $r=1$ for $1 \le r \le 2$ and $\alpha_2(1;1)$ has the asymptotic behavior (\ref{alpha-asym}), 
the product $2^d \alpha_2(1;1) \phi$ vanishes no more
slowly than $(6/\pi)^{1/2}/[(4/3)^{d/2}d^{1/2}]$ in the limit $d \rightarrow \infty$ for
$0 \le \phi \le 1/2^d$, and therefore  $g_2(r)$ tends  to 
$\Theta(r-1)$ exponentially fast. In summary, we see again that spatial correlations
that exist in low dimensions for $r >1$ completely vanish in the limit
$d \rightarrow \infty$. Moreover,  this is yet another disordered
packing model in which the  high-dimensional
asymptotic behavior corresponds to the low-density asymptotic behavior.

The corresponding $n$-particle correlation function $g_n$, defined by (\ref{nbody}),
in the low-density limit \cite{Sa58} is given by
\begin{equation}
g_n({\bf r}_{12}, \ldots, {\bf r}_{1n})= \prod_{i<j}^n g_2(r_{ij})
\Big[1+ 2^d \alpha_n({\bf r}_{12},\ldots, {\bf r}_{1n};1) \phi + {\cal O}(\phi^2)\Big]
\end{equation}
where $ \alpha_n({\bf r}_{12},\ldots, {\bf r}_{1n};R)$ is the intersection volume
of $n$ congruent spheres of radius $R$ (whose centers are located at
${\bf r}_1,\ldots, {\bf r}_n$, where ${\bf r}_{ij}={\bf r}_j -{\bf r}_i$
for all $1 \le i <j \le n$) divided by the volume of a sphere of
radius $R$. The scaled intersection volume $\alpha_n({\bf r}_{12},\ldots, {\bf r}_{1n};R)/n$
has the range $[0,1]$. Now since  $\alpha_2(r_{ij},1) \ge \alpha_n({\bf r}_{12},\ldots, {\bf r}_{1n};1)$
for any pair distance $r_{ij}=|{\bf r}_{ij}|$ such that $1 \le i <j \le n$, 
then it follows from the analysis above that in the limit $d \rightarrow \infty$ for
$0 \le \phi \le 1/2^d$
\begin{equation}
g_2(r) \sim \Theta(r-1)
\label{g2-eq-asym}
\end{equation}
and
\begin{equation}
g_n({\bf r}_{12}, \ldots, {\bf r}_{1n}) \sim \prod_{i<j}^n g_2(r_{ij}).
\label{gn-eq-asym}
\end{equation}
Again, as in the generalized RSA example with $\kappa=1$, $g_n$ factorizes into 
products involving only $g_2$'s in the limit $d\rightarrow \infty$.
Moreover, we should also note that the standard RSA process
(generalized RSA process with $\kappa=0$) has precisely the same asymptotic
low-density behavior as the standard Gibbs hard-sphere model \cite{To02a}. More precisely,
these two models share the same low-density expansions of the $g_n$
through terms of order $\phi$ and therefore the same asymptotic expressions
(\ref{g2-eq-asym}) and (\ref{gn-eq-asym}).

\subsection{Decorrelation Principle}

 The previous two examples illustrate two important and related
asymptotic properties that are expected to apply to all disordered packings: 
\begin{enumerate}
\item  the high-dimensional asymptotic behavior of $g_2$ is the same as the asymptotic
behavior in the low-density limit for any finite $d$, i.e., {\it unconstrained} spatial
correlations, which exist for positive densities at fixed $d$, vanish
asymptotically for pair distances beyond the hard-core diameter in the high-dimensional limit; 
\item and $g_n$ for $n \ge 3$ asymptotically can be inferred from
a knowledge of only the pair correlation function $g_2$ and number density $\rho$.
\end{enumerate}
What is the explanation for these two related asymptotic properties? Because we know from  
the Kabatiansky and Levenshtein (1979) asymptotic upper bound on the maximal density
that $\phi$ must go to zero at least as fast as $2^{-0.5990d}$ for large $d$,
unconstrained spatial correlations between spheres must vanish, i.e., statistical
independence is established. (An example
of {\it constrained} spatial correlations is described below.) Such a decorrelation
means that the $g_n$ for $n \ge 3$ are determined entirely from a knowledge of the decorrelated pair
correlation function $g_2$. In the specific examples that we considered, the $g_n$
factorize into products involving only $g_2$'s, but there may be other decompositions.
For example, the $g_n$ for $n \ge 3$ can be functionals that only involve
$\rho$ and $g_2$. We will call the two asymptotic properties the {\it decorrelation principle}
for disordered packings. This principle as well as other results described
in Section 5 leads us to a conjecture concerning
the existence of disordered sphere packings in high dimensions,
which we state in Section 5.1.

\begin{figure}[bt]
\centerline{\psfig{file=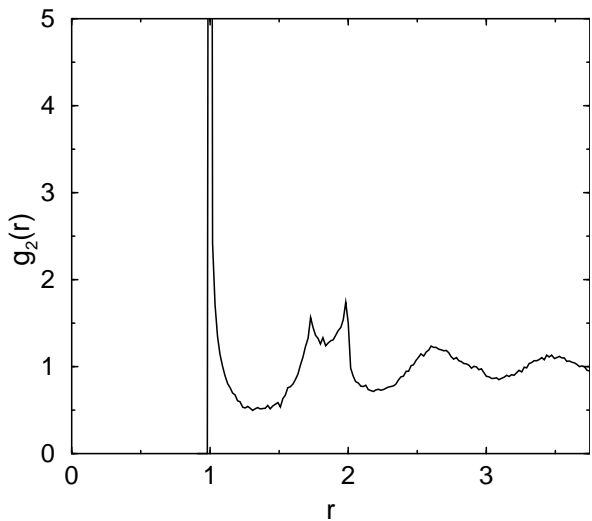,height=2.8in,clip=}}
\caption{The radial distribution function
of  a three-dimensional  packing of spheres near
the maximally random jammed state \cite{To00b,To02a} at a density $\phi=0.64$
as obtained from computer simulations \protect\cite{To02c}.
The delta function contribution at $r=1$ (of course, not explicitly shown)
corresponds to an average kissing number of about six.}
\label{MRJ}
\end{figure}

An example of constrained spatial correlations that would not vanish asymptotically
is illustrated in Fig. \ref{MRJ}, where we show the pair correlation
function $g_2(r)$ for a three-dimensional sphere packing near the so-called {\it maximally random jammed
state} \cite{To00b,To02a}. A special feature of $g_2(r)$ for a maximally random jammed
packing is a delta-function contribution at $r=1$, which reflects the fact that the 
{\it average kissing number} (i.e., average number of contacting particles per particle)
is effectively six for this collectively jammed packing, meaning that
the packing is {\it isostatic} \cite{Do05c}. A positive average kissing number is required if the packing
is constrained to be jammed and in $\mathbb{R}^d$ this means that the average
kissing number is $2d$ for either collective or strict jamming \cite{Do05c}. Isostatic packings are jammed
packings with the minimum number of contacts for a particular jamming category.
According to the Decorrelation Principle, as $d$ tends to infinity, 
$g_2$ for a maximally random jammed packing would retain this delta-function contribution but
the unconstrained spatial correlations beyond $r=1$ would vanish. Of course, the
manner in which the $g_2$ shown in Fig. \ref{MRJ} approaches the asymptotic
limit of a step function $\Theta(r-1)$ plus a delta-function contribution
at $r=1$ is crucial. We note that maximally random jammed packings 
contain about 2--3\% {\it rattler} spheres, i.e., spheres
trapped in a cage of jammed neighbors but free to move within the
cage.

%\begin{postulate}
%Suppose there exists a translationally invariant point process
%characterized by the  $n$-particle correlation functions $g_n({\bf r}_{12},\ldots,{\bf r}_{1n})$
%that belong to the class of tempered distributions in $\mathbb{R}^d$ for any $n \ge 3$
%and are entirely determined by the number density $\rho$ and pair correlation
%function $g_2({\bf r})$. The density range over which such a point process
%exists is the one in which $\rho$ and $g_2({\bf r})$ satisfy the standard two
%nonnegativity conditions (\ref{cond1}) and (\ref{cond2}).
%\label{axiom}
%\end{postulate}

\section{New Approach to Lower Bounds}
\label{new}

The salient ideas behind our new approach to the derivation of lower bounds on $\phi_{\mbox{\scriptsize max}}$
were actually laid out in our earlier work \cite{To02c}. The main objective
of that work was to study sphere packings in three dimensions in which
long-range order was suppressed and short-range order was controlled (i.e.,
disordered sphere packings in $\Re^3$) using so-called $g_2$-invariant processes.
A {\it $g_2$-invariant process} is one in which a given nonnegative pair correlation $g_2({\bf r})$
function remains invariant for all ${\bf r}$ over the range of
densities
\begin{equation}
0 \le \phi \le \phi_*.
\end{equation}
The terminal density $\phi_*$ is the maximum achievable density
for the $g_2$-invariant process subject to satisfaction of the structure factor 
$S(\bf k)$  inequality (\ref{cond2}). A five-parameter test family of $g_2$'s had been considered,
which incorporated the known features of core exclusion, contact pairs, and damped oscillatory short-range
order beyond contact.  The problem of finding the maximal packing fraction $\phi_*$ was posed
as an optimization problem: maximize $\phi$ over the set of parameters subject
to the constraints (\ref{cond1}) and (\ref{cond2}). We noted in passing
that when the damped-oscillatory contribution to $g_2$ was set equal to zero, 
the optimization problem could be solved analytically for all space dimensions,
leading to a terminal density $\phi_*=(d+2)/2^{d+1}$. Under the assumption
that such a $g_2$ was a realizable packing, we also observed that this $\phi_*$
was a lower bound on the maximal density for 
any sphere packing [i.e., $\phi_{\mbox{\scriptsize max}} \ge (d+2)/2^{d+1}$] because the
terminal density would have been higher by including the damped-oscillatory contribution to $g_2$.
This conjectural lower bound was noted to provide linear improvement
over Minkowski's lower bound, but we were not aware of Ball's similar lower bound \cite{Ball92}
at the time. Since our original 2002 paper, we also learned about other necessary conditions
for the realizability of a point process for a given number density $\rho$ and
$g_2$, such as Yamada's condition (\ref{yamada}). In any event, our brief remarks about
lower bounds on sphere packings were not intended to be mathematically rigorous.

It is our intent here to make our optimization methodology
to obtain lower bounds on $\phi_{\mbox{\scriptsize max}}$ more mathematically precise, especially
in light of recent developments and the considerations of the previous two sections.
We then apply the optimization procedure to provide alternative
derivations of previous lower bounds as well as a new bound.

We will consider those ``test" $g_2(r)$'s that are distributions on $\mathbb{R}^d$ depending
only on the radial distance $r$ such that $h(r)=g_2(r)-1$. 
For any test $g_2(r)$, we want to maximize 
the corresponding density $\phi$ satisfying the following three conditions:

\noindent (i) \hspace{0.25in} $g_2(r)  \ge 0 \qquad \mbox{for all}\quad  r,$

\noindent (ii) \hspace{0.2in} $\mbox{supp}(g_2) \subseteq \{r: r \ge 1\}$,

\noindent (iii)
\begin{displaymath}
\hspace{-0.25in}S(k)= 1+ \rho\left(2\pi\right)^{\frac{d}{2}}\int_{0}^{\infty}r^{d-1}[g_2(r)-1]
\frac{J_{\left(d/2\right)-1}\!\left(kr\right)}{\left(kr\right)^{\left(d/2\right
)-1}} \, dr \ge 0  \qquad \mbox{for all}\quad  k.
\end{displaymath}
We will call the maximum density the {\it terminal density} and denote it by $\phi_*$.

\noindent{\bf Remarks:}
\smallskip 

\noindent 1. The conditions (i)--(iii) are just recapitulations of (\ref{support}), (\ref{cond1}),
and (\ref{cond2}) for this class of test functions. We will call
condition (ii) the {\it hard-core} constraint.
\smallskip

\noindent 2. When there exist sphere packings with $g_2$ satisfying conditions
(i)--(iii) for $\phi$ in the interval $[0,\phi_*]$, then we have the lower
bound on the maximal density given by
\begin{equation}
\phi_{\mbox{\scriptsize max}} \ge \phi_*.
\label{true-bound-phi}
\end{equation}
The best lower bound would be obtained if one could probe the entire class of
test functions. In practice, we will consider here only a small subset
of test functions and in particular those that are amenable to exact asymptotic
analysis. In some instances, we will associate with the terminal density $\phi_*$ an
optimized {\it average kissing number} $Z_*$. Thus, whenever inequality (\ref{true-bound-phi})
applies, the {\it maximal kissing number} $Z_{\mbox{\scriptsize max}}$ is bounded from below
by $Z_*$, i.e.,
\begin{equation}
Z_{\mbox{\scriptsize max}} \ge Z_*.
\label{true-bound-Z}
\end{equation}
In the next subsection, we put forth a conjecture that states when
the conditions (i)--(iii) are necessary and sufficient for
the existence of disordered sphere packings.

\smallskip

\noindent 3. Remarkably, the optimization problem defined above is identical
to one formulated by Cohn \cite{Co02}. In particular, it is the dual
of the primal infinite-dimensional linear program that Cohn employed with Elkies \cite{Co03}
to obtain upper bounds on the maximal packing density. One need 
only replace $S(k)$ with ${\hat g}- c\delta(k)$, where 
$c$ plays the role of number density,
$g$ is a tempered distribution and ${\hat g}$ is its Fourier
transform in Cohn's notation. Thus, even if  there does 
not exist a sphere packing with $g_2$ satisfying conditions
(i)--(iii), our formulation has implications for upper bounds
on $\phi_{\mbox{\scriptsize max}}$, which we discuss in
Section \ref{discuss}. For finite-dimensional linear programs (and many infinite-dimensional
ones) there is no ``duality gap," i.e., the optima of the primal
and dual programs are equal. However, in this infinite-dimensional
setting, it is not clear how to prove that there is no duality gap \cite{Co02}.
Therefore, it is rigorously true that the terminal density $\phi_*$ can never exceed the
Cohn-Elkies upper bound, which is a desirable feature of our formulation,
otherwise the terminal density could never correspond to a rigorous lower bound.

\bigskip

We will show that for the test radial distribution functions
considered in this paper, the Yamada condition (\ref{yamada}) is only relevant
in one dimension, and even then in just some cases. A remark about the Yamada condition  for sphere
packings is in order here. In earlier work \cite{To03a}, we observed that for any
sphere packing of congruent spheres, the number variance
for a spherical window of radius $R$ defined by (\ref{variance})  obeys the lower bound
\begin{equation}
\sigma^2(R) \ge 2^d \phi R^d\Big[1-2^d\phi R^d\Big]
\label{lower}
\end{equation}
for any $R$. This is a tight bound for sufficiently small $R$
and is exact for $R \le 1/2$. However, we note here that provided that 
$2^d \phi R^d \le 1$, the Yamada lower bound (\ref{yamada}) and lower bound
(\ref{lower}) are identical. Thus, the Yamada lower bound
for any sphere packing only needs to be checked for $R>R_0$, where
\begin{equation}
R_0= \frac{1}{2 \phi^{1/d}}.
\label{yamada-int}
\end{equation}

\subsection{Existence of Disordered Packings in High Dimensions}

We have seen that a necessary condition for the existence of a 
translationally invariant point process with a specified positive $\rho$ 
and nonnegative $g_2$ is that $S({\bf k})$ is nonnegative [cf. (\ref{cond2})].
In other words, given $\rho$ and $g_2$, it does not mean that there
are some higher-order functions $g_3, g_4, \ldots$ for which these one- and two-point
correlation functions hold. The function $g_2$ specifies how frequently
pair distances of a given length occur statistically in $\mathbb{R}^d$.  The third-order
function $g_3$ reveals how these pair separations are linked into triangles.  This
additional information generally cannot be inferred from the
knowledge of $\rho$ and $g_2$ alone, however. The fourth-order
function $g_4$ controls the assembly of triangles into tetrahedra (and is the
lowest-order correlation function that is sensitive to chirality)
but  $g_4$ cannot be determined by only knowing $\rho$, $g_2$ and $g_3$.
In general, $g_n$ for any $n \ge 3$ is not completely determined
from a knowledge of the lower-order correlations functions alone.
This is to be contrasted with general stochastic processes in which
nonnegativity of first- and second-order statistics (mean and autocovariance)
are necessary and sufficient to establish existence because one can always find
a Gaussian process with such given first- and second-order statistics.
For a Gaussian process, first- and second-order statistics determine
all of the high-order statistics.

There are a number of results that suggest it is reasonable to conclude that 
the generally  necessary nonnegativity conditions for the existence of a 
disordered sphere packing become necessary and sufficient for sufficiently large $d$.
First, the decorrelation principle of the previous section states that unconstrained
correlations in disordered sphere packings vanish asymptotically in high dimensions
and that the $g_n$ for any $n \ge 3$ can be inferred entirely from a knowledge
of $\rho$ and $g_2$.  Second, as we noted in Section \ref{disordered}, the necessary 
Yamada condition appears to only have relevance in very low dimensions.
Third, we will demonstrate below that other new necessary conditions
also seem to be germane only in very low dimensions. Fourth, we will
describe  numerical constructions of configurations of disordered sphere
packings on the torus corresponding to certain test radial distributions functions
in low dimensions for densities up to the terminal density. Finally, we will show that
certain test radial distributions functions recover the asymptotic forms of known
rigorous bounds. In light of these results, we propose the following conjecture:

\begin{conjecture}
 For sufficiently large $d$, a hard-core nonnegative tempered distribution $g_2({\bf r})$ 
 that satisfies $g_2({\bf r}) = 1 + {\cal O}(|{\bf r}|^{-d-\varepsilon})$ for some
 $\varepsilon>0$ is a pair correlation function
 of a translationally invariant disordered sphere packing in $\mathbb{R}^d$ at number density $\rho$ 
 if and only if $S({\bf k}) \equiv 1+\rho {\tilde h({\bf k})} \ge 0$. The maximum achievable density
 is the terminal density $\phi_*$.
\label{conj}
\end{conjecture}

\noindent{\bf Remarks:}
\smallskip 

\noindent 1. A weaker form of this conjecture would replace the phrase 
``for sufficiently large $d$" with ``in the limit $d \rightarrow \infty$."

\noindent 2. Employing the aforementioned optimization procedure with a certain
test function $g_2$ and this conjecture, we obtain in what follows
conjectural lower bounds that yield the long-sought
asymptotic exponential improvement on Minkowski's bound.
Before obtaining this result, we first apply the procedure
to two simpler test functions that we examined in the past.

\subsection{Step Function}

The simplest possible choice for a radial distribution function
corresponding to a disordered packing is the following unit step function: 
\begin{equation}
g_2(r)=\Theta(r-1).
\label{step1}
\end{equation}
This states that all pair distances beyond the hard-core diameter
are equally probable, i.e., spatial correlations vanish identically.
The corresponding structure factor [cf. condition (iii)] for this test function
in any dimension $d$ is given by \cite{To02c}
\begin{equation}
S(k)=1 -\frac{\phi 2^{3\nu} \Gamma(1+\nu) }{k^{\nu}} J_{\nu}(k),
\end{equation}   
where $\nu=d/2$. Since there are no parameters to be optimized here, the terminal density $\phi_*$ is readily obtained by determining
the highest density for which the condition (\ref{cond2}) is satisfied, yielding
\begin{equation}
\phi_*=\frac{1}{2^d}.
\label{term1}
\end{equation}

Now we show that the Yamada condition (\ref{yamada}) is satisfied in any dimension
for $0 \le \phi \le 2^{-d}$.
Consider the more general class of radial distribution functions:
\begin{equation}
0 \le g_2(r) \le 1 \qquad \mbox{for}\quad r>1. 
\label{h-class}
\end{equation}
The test function (\ref{step1})  belongs to this class. Note that for any dimension, the scaled intersection
volume given by (\ref{series}) obeys the inequality
\begin{equation}
\alpha_2(r;R) \le 1 -\frac{r}{2R} \qquad {\mbox{for}}\quad  0 \le r \le 2R,
\label{alpha-class}
\end{equation}
where the equality applies when $d=1$. For $g_2(r)$ satisfying (\ref{h-class}),
relation (\ref{variance}) and inequality (\ref{alpha-class}) yield
the following lower bound for any $d$:
\begin{equation}
\sigma^2(R) \ge 2^d \phi R^d \Bigg[1+d2^d \phi\int_{0}^{2R} r^{d-1} h(r)\Big[1-\frac{r}{2R}\Big]
\, dr\Bigg].
\label{var}
\end{equation}

At $\phi=1/2^d$, the lower bound (\ref{var}) for the test function (\ref{step1}) is given by
\begin{equation}
\sigma^2(R) \ge \frac{d}{2(d+1)} R^{d-1}
\end{equation}
and because $R_0=1$ [cf. (\ref{yamada-int})], we only need to consider $R>1$. In particular,
the right side of this inequality is smallest at $R=1$ so that
\begin{equation}
\sigma^2(R) \ge \frac{d}{2(d+1)}.
\end{equation}
Since $\sigma^2(R) \ge 1/4$ for $d \ge 1$, Yamada's condition
is satisfied for all $R$ for the step function (\ref{step1}) at $\phi=1/2^d$
as well as all $\phi <1/2^d$. 

We already established in Section \ref{disordered} that there exist sphere packings that 
asymptotically have radial distribution functions given by the simple
unit step function (\ref{step1}) for $\phi \le 2^{-d}$. Nonetheless,
invoking Conjecture \ref{conj} and terminal density specified by (\ref{term1})
implies the asymptotic lower bound on the maximal density is given by
\begin{equation}
\phi_{\mbox{\scriptsize max}} \ge \frac{1}{2^d},
\end{equation}
which provides an alternate derivation of the elementary bound (\ref{sat}).

Using numerical simulations with a finite but large
number of spheres on the torus, we have been able to construct particle configurations
in which the radial distribution function (sampled at discretized pair distances)
is given by the test function (\ref{step1}) in one, two and three dimensions for densities
up to the terminal density \cite{Cr03,Ou06}. The existence of such
a discrete approximation to (\ref{step1}) of course is not conclusive proof of the 
existence of such packings in low dimensions, but they
are suggestive that the standard nonnegativity
conditions may be sufficient to establish existence in this case 
for densities up to $\phi_*$.

\subsection{Step Plus Delta Function}

An important feature of any dense packing is that the particles
form contacts with one another. Ideally, one would like
to enforce strict jamming (see Introduction). The probability
that a pair of particles form such contacts at the pair
distance $r=1$ for the test function (\ref{step1}) is strictly zero. 
Accordingly, let us now consider the test radial distribution function given by
the previous test function plus a delta-function contribution
as follows:
\begin{equation}
g_2(r)=\Theta(r-1)+ \frac{Z}{s_1(1)\rho}\delta(r-1).
\label{step2}
\end{equation}
Here $s_1(r)$ is the surface area of a  $d$-dimensional sphere of radius $r$ given by (\ref{surf})
and $Z$ is a parameter, which is the average kissing number.  Because we allow for interparticle contacts via
the second term in (\ref{step2}), the terminal density is expected to be greater
than $2^{-d}$, which will indeed be the case.
The corresponding structure factor [cf. (iii)] for this test function
in any dimension $d$ is given by \cite{To02c}
\begin{equation}
S(k)=1 -\frac{\phi 2^{3\nu} \Gamma(1+\nu) }{k^{\nu}} J_{\nu}(k)+\frac{Z 2^{\nu}\Gamma(1+\nu)
}{d k^{\nu-1}} J_{\nu-1}(k),
\end{equation}      
where $\nu=d/2$. 
The structure factor for small $k$ can be expanded in a MacLaurin series
as follows:
\begin{equation}
S(k)=1+(Z-2^d\phi)+\left[\frac{2^{d-2}\phi}{1+d/2}-\frac{Z}{2d}\right] k^2+
{\cal O}(k^4).
\end{equation}
The last term changes sign if $Z$ increases past $2^d \phi d/(d+2)$.
At this crossover point,
\begin{equation}
S(k)=1-\frac{2^{d+1}}{d+2} \phi+ {\cal O}(k^4)
\end{equation}
Under the constraint that the minimum of $S(k)$ occurs at $k=0$, then we have the exact results
\begin{equation}
\phi_*=\frac{d+2}{2^{d+1}}, \qquad Z_*=\frac{d}{2}.
\label{term2}
\end{equation}
We see that at the terminal density, the average kissing number $Z_*={d}/{2}$, which does not
even meet the local jamming criterion described in the Introduction. 

The Yamada condition (\ref{yamada}) is violated only for $d=1$ for the test function
(\ref{step2}) at the terminal density specified by (\ref{term2}). It is easy
to verify directly that the Yamada condition becomes less restrictive as the
dimension increases from $d=2$.
Interestingly, we have also shown via numerical simulations
that there exist sphere packings possessing radial distribution functions 
given by the test function (\ref{step2}) (in the discrete approximation)
in two and three dimensions for densities
up to the terminal density \cite{Ou06}. This is suggestive that Conjecture
\ref{conj} for this test function may in fact be stronger than is required.

In the high-dimensional limit, we invoke  Conjecture \ref{conj} and the terminal density given 
by (\ref{term2}), yielding the conjectural lower bound
\begin{equation}
\phi_{\mbox{\scriptsize max}} \ge \frac{d+2}{2^{d+1}}.
\label{lower-2}
\end{equation}
This lower bound provides the same type of linear improvement over Minkowski's lower
bound as does Ball's lower bound \cite{Ball92}.

\subsection{Step Plus Delta Function with a Gap}

 The previous test function (\ref{step2}) provided an optimal average kissing
number $Z_*=d/2$ that did not even meet the local jamming criterion. Experience
with disordered jammed packings in low dimensions reveals that the kissing
number as well as the density can be substantially increased if there
is there is a low probability of finding noncontacting particles from a typical 
particle at radial distances just larger than the nearest-neighbor distance.
This small-distance negative correlation is clearly manifested
in the graph of $g_2(r)$ for the three-dimensional maximally random jammed
packing (Figure \ref{MRJ}) for values of $r$ approximately between
1.1 and 1.5. We would like to idealize this small-distance negative correlation
in such a way that it is amenable to exact asymptotic analysis.
Accordingly,  we consider a test radial distribution function that is similar to
the previous one [cf. (\ref{step2})] but one in which there is
a gap between the location of the unit step function and the
delta function at finite $d$, i.e.,
\begin{equation}
g_2(r)=\Theta(r-\sigma)+\frac{Z}{s_1(1)\rho}\delta(r-1).
\label{step3}
\end{equation}
The expression contains two adjustable parameters, $\sigma \ge 1$ and $Z$,
which must obviously be constrained to
be nonnegative.  According to the decorrelation principle
of Section \ref{disordered}, the location of the step function $r=\sigma$
must approach unity asymptotically, i.e., it must approach the previous test function (\ref{step2}).
However, as we have emphasized, the manner
in which the test function (\ref{step3}) approaches (\ref{step2})
is crucial. Indeed, we will see that the presence of 
a gap between the unit step function and delta function will indeed lead
asymptotically to substantially higher terminal densities.

The structure factor is given by
\begin{equation}
S(k)=1 -\frac{ 2^{3\nu}\phi\sigma^{d} \Gamma(1+\nu) }{(k\sigma)^{\nu}} J_{\nu}(k\sigma)
+  \frac{2^{\nu}Z \Gamma(1+\nu) }{d k^{\nu-1}} J_{\nu-1}(k).
\label{S}
\end{equation}               
The goal now is to find the optimal values of the 
the adjustable nonnegative parameters $Z$ and $\sigma$ 
that maximize the density $\phi$ subject to the constraint
 (iii).  This search in two-dimensional parameter space can be reduced
by imposing the further condition that a minimum of the structure factor
occurs at $k=0$. The MacLaurin expansion of expression (\ref{S}) gives
\begin{equation}
S=1+[Z-(2\sigma)^d \phi]+\Bigg[ \frac{2^{d-1}\sigma^{d+2}\phi}{d+2}-\frac{Z}{2d}\Bigg] k^2+
{\cal O}(k^4).
\end{equation}
Requiring that a zero of $S(k)$ occurs at the origin (hyperuniformity) such that
the quadratic coefficient is nonnegative implies the restrictions
\begin{equation}
Z=(2\sigma)^d\phi -1
\label{Z}
\end{equation}
and 
\begin{equation}
(2\sigma)^d\phi[d\sigma^2-(d+2)]+ d+2 \ge 0.
\end{equation}
Combination of (\ref{S}) and (\ref{Z}) yields the structure factor as
\begin{equation}
S(k)=1-c_1(d) \frac{J_{\nu}(k\sigma)}{(k\sigma)^{\nu}}
+c_2(d) \frac{J_{\nu-1}(k)}{(k\sigma)^{\nu-1}},
\label{S3}
\end{equation}
where  the $d$-dependent coefficients $c_1(d)$ and $c_2(d)$ are given by
\begin{equation}
c_1(d)=\phi\sigma^d 2^{3\nu}\Gamma(1+\nu)
\label{c1}
\end{equation}
\begin{equation}
c_2(d)=[(2\sigma)^d\phi -1]2^{\nu}\frac{\Gamma(1+\nu)}{d}.
\label{c2}
\end{equation}
Now the  problem reduces to finding the optimal value of the parameter
$\sigma(d)$ as a function of the space dimension $d$ that maximizes the density $\phi$ subject to 
the constraint (iii).  It will be shown below that the optimal $\sigma$ is of order unity and
approaches unity in the limit $d \rightarrow \infty$. It immediately follows
from (\ref{S3}) and the asymptotic properties of the Bessel functions of fixed order that
 $S(k) \rightarrow 1$ for $k \rightarrow \infty$. 

In general, $S(k)$ will possess multiple minima and thus we want to ensure that the values of $S(k)$ 
at each of these minima are all nonnegative. To find the minima of $S(k)$, we set its first derivative 
to zero, yielding the relation
\begin{equation}
\frac{c_1(d)}{\sigma^{\nu-1}} J_{\nu+1}(k\sigma)=c_2(d) k J_{\nu}(k),
\label{minima1}
\end{equation}
where we have used the identity 
\begin{equation}
\frac{d}{dx} \Bigg[ \frac{J_{\nu}(x)}{x^{\nu}}\Bigg]=-\frac{J_{\nu+1}(x)}{x^{\nu}}.
\end{equation}

\begin{figure}
\centerline{\psfig{file=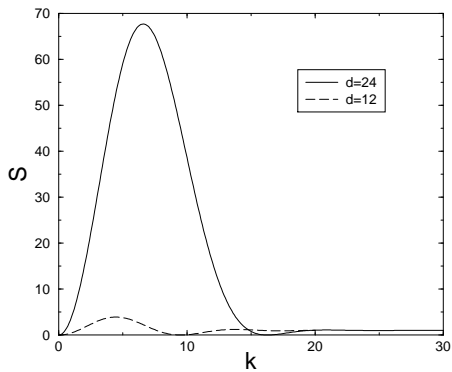,height=2.0in}}
\caption{The optimized structure factor for $d=12$ and $d=24$.}
\label{12-24}
\end{figure}

For sufficiently small $d$ ($d \le 200$), the  search procedure is carried out numerically using
Maple, and is made more efficient by exploiting
the fact that the minima of $S(k)$ occur at the real solutions  of (\ref{minima1}).
Figure \ref{12-24} shows the optimized structure factor for $d=12$ and $d=24$.
Our numerical examination of $S(k)$ for a wide range of $d$ values has consistently shown
that the first minimum for positive $k$ is the deepest one. Although we have not proven this rigorously,
we assume that this is a general result.

We should note that the Yamada condition (\ref{yamada}) is 
violated only for $d=1$ for the test function
(\ref{step3}) for the terminal density $\phi_*$ and associated optimized 
parameters $\sigma_*$ and $Z_*$ [calculated via (\ref{Z})]. One can
again verify directly that the Yamada condition becomes less restrictive as the
dimension increases from $d=2$. However, although the test function
(\ref{step3}) for $d=2$ with optimized parameters $\phi_*=0.74803$,
$\sigma_*=1.2946$ and $Z_*=4.0148$ satisfies the Yamada condition,
it cannot correspond to a sphere packing because it violates
local geometric constraints specified by $\sigma_*$ and $Z_*$. 
Specifically, for an average kissing number of 4.0148, there must be particles 
that are in contact with at least five others. 
But no arrangement of the five exists which is consistent with the assumed
pair correlation function (step plus delta function with a gap from 1 to
1.2946).  Simple geometric considerations show that either some pairs of the
five would be forced into the gap, or they would be restricted to fixed
separations that would correspond to undesired delta functions beyond the
gap.  To our knowledge, this is the first example of a test radial
distribution function that satisfies the two standard non-negativity
conditions (\ref{cond1}) and (\ref{cond2}) and the Yamada condition (\ref{yamada}), but cannot
correspond to a point process.  Thus, there is at least one previously
unarticulated necessary condition that has been violated in the low
dimension $d=2$.

\begin{table}[bt]
\centering
\caption{Optimized parameters $\sigma_*$, $Z_*$ and $\phi_*$,
and the ratio  $\frac{2^{d+1}\phi_*}{d+2}$, which is the relative improvement
of the terminal density over the gapless-test-function terminal
density [cf. (\ref{term2})].\label{parameters}} 
\footnotesize
\begin{tabular}{c|c|c|c|c}
\multicolumn{5}{c}{~} \\ \hline\hline
d  &     $\sigma_*$&     $Z_*$     &      $\phi_*$ &           $\frac{2^{d+1}\phi_*}{d+2}$  \\ \hline
3  &     1.246997&    7.932582  &   0.5758254    &  1.842641  \\ %
4  &     1.212589&   13.71016  &    0.4252472    &  2.267985  \\ %
5  &     1.186929&   21.97918   &   0.3048322    &  2.787037 \\ %
6  &     1.167000&   33.53884&      0.2136444    &  3.418310 \\ %
7  &     1.151106&   49.42513&      0.1471058    &  4.184343  \\%
8  &      1.137967&  70.88348&     0.09985085&      5.112364  \\ %
%9  &       1.12691&    99.52347&     0.0669876 &        6.24\\
%11 &      1.11016 &   188.4056&     0.0292969 &        9.23\\
%12 &      1.1031  &   253.6286  &      0.01915   &       11.21\\
%13 &      1.098   &   341.48539041&     0.0124    &       13.54\\
%15 &      1.087   &   589.93919726&     0.00516   &       19.89\\
%17 &      1.0787  &   999.16940382&     0.002105  &       29.04\\
%19 &      1.072   &  1659.07438461&     0.000845  &       42.19\\
%21 &      1.0663  &  2710.30387085&     0.0003358 &       61.23\\
%23 &      1.062   &  4382.54108210&     0.000131  &       87.91\\
24 &      1.058992&  5473.546 &     8.245251e-05&   106.4095   \\ %
%25 &      1.0572  &  6927.54228931&     0.0000514 &      127.76\\
%27 &      1.0535  & 10891.82003510&     0.00001987&      183.92\\
%29 &      1.05037 &               &     0.00000763&      264.27\\
%31 &      1.0479  & 26651.49715   &     0.00000291&      378.74\\
%33 &      1.04642 &               &     0.000001095&     537.48\\
%35 &      1.0442  & 64479.59117   &     0.000000413&     767.06\\
36  &     1.041611 & 76521.15    &    2.566299e-07&  928.1828   \\%
%37 &      1.0422  &               &     0.000000155&    1092.46\\
%39 &      1.0393  &               &     5.89  e-08  &   1579.54\\
%41 &      1.0375  &               &     2.20  e-08  &   2250.16\\
%43 &      1.0359  &               &     8.20  e-09  &  3205.69 \\
%45 &      1.0342  & 486883.8528   &     3.047 e-09  & 4561.99 \\
%47 &      1.0330  &               &     1.128 e-09  & 6479.67\\
%49 &      1.0317  &               &     4.166 e-10  &   9197.06\\
%51 &      1.0306  & 1.607810769e06&     1.535 e-10 &   13043.44\\
%53 &      1.0295  &               &     5.644 e-11 &   18486.05\\
%55 &      1.0286  &               &     2.070 e-11 &   26168.28\\
56 &      1.028036& 4.248007e06&    1.253255e-11&  31140.19   \\ %
%57 &      1.0276  &               &     7.581 e-12 &   37035.16     \\
%59 &      1.0268  &               &     2.770 e-12 &   52353.98     \\
60 &       1.026330& 9.179315e06&    1.674130e-12&  62262.60  \\%
%61 &      1.0261  & 1.121300497e07&     1.010 e-12 &   73933.38     \\
%63 &      1.0252  &               &     3.682 e-13 &  104493.71     \\
64 &      1.024823&  1.968233e07 &    2.221414e-13& 124175.32    \\%
%65 &      1.0246 &                &    1.339 e-13  &  147463.82     \\
%67 &      1.0238 &                &    4.864 e-14  &  208057.88     \\
%69 &      1.0232 &                &    1.764 e-14  &  293318.82     \\
%71 &      1.0226 &                &    6.389 e-15  &  413304.10     \\
80 &      1.020211& 3.908042e08 &    6.521679e-17& 1.922982e06 \\  %
100&      1.016421& 1.478804e10 &    2.288485e-21& 5.688234e08 \\    
125&      1.013311& 1.246172e12 &    5.610270e-27& 3.758024e09 \\    
150&      1.011214& 9.698081e13 &    1.275632e-32& 2.319290e11 \\
175&      1.009671& 7.086019e15 &    2.745830e-38& 1.485866e13 \\
200&      1.008510& 4.959086e17 &    5.667098e-44& 9.016510e14 \\   \hline\hline %
\end{tabular}
\end{table}

In three dimensions one obtains $\phi_*=0.5758254$, $\sigma_*=1.246997$, and
$Z_*=7.932582$.  The last of these requires that some nonzero fraction of the
spheres must have at least 8 contacting neighbors.  We have verified that
valid arrangements of both 8 and 9 contacts are possible, thereby avoiding
the analog of  the violation encountered in $d=2$.  
As is the case with the Yamada condition (\ref{yamada}), 
this additional necessary condition appears to lose
relevance as $d$ increases. 

The terminal density $\phi_*$ and the associated
optimized parameters $\sigma_*$ and $Z_*$ are listed
in Table \ref{parameters} for selected values of the space dimension between $d=3$ and $d=200$.
Note that for $d \le 56$, the terminal density lies below the density of the densest
known packing. For $d=56$, the densest arrangement is a lattice (designated
by $L_{56,2}(M)$ \cite{Ne98}) with density $\phi=2.327670 \times 10^{-11}$,
which is about twice as large as $\phi_*$, as shown in the table. However,
for $d >56$,  $\phi_*$ can be larger than the density of the densest
known arrangement. For $d=60$,  the densest known packing is again a lattice (designated
by $L_{56,2}(M)$ \cite{Co93}) with density $\phi=2.966747 \times 10^{-13}$,
which is about five times smaller than $\phi_*$, as shown in the table.
The next dimension for which data is available is $d=64$, where
the densest known packing is the $Ne_{64}$ lattice \cite{Ne98}
with density $\phi=1.326615 \times 10^{-12}$, which is about six times larger than 
$\phi_*$. The table also reveals exponential improvement of the terminal density
$\phi_*$ over the one for the gapless case, i.e., $\phi_*=(d+2)/2^{d+1}$.
The crucial question is whether such exponential improvement persists
in the high-dimensional limit.

To obtain an asymptotic expression for $\phi_*$ for large $d$, we use the fact
that $(2\sigma)^d \phi \gg 1$, implying that $c_1(d)/c_2(d) \rightarrow d$
[cf. ({\ref{c1}) and (\ref{c2})]. Therefore,
the minima of $S(k)$ for large $d$ are the solutions of
\begin{equation}
\frac{J_{\nu+1}(k\sigma)}{\sigma^{\nu-1}}= \frac{k}{d} J_{\nu}(k).
\label{minima2}
\end{equation}
We see that the locations of the minima depend only on $\sigma$ (not on $\phi$).
The deepest minimum of $S(k)$, after the one at $k=0$, is a zero 
and occurs at the  wavenumber $k=k_{min}$. [This characteristic is
true in any dimension (see Figure \ref{12-24}).]
Therefore, $S(k=k_{min})=0$,  $c_2(d) =c_1(d)/d$, and
relation (\ref{S}) gives the condition
\begin{equation}
\frac{c_1(d)}{k_{min}^{\nu}} \Delta_{\nu}(k_{min})=
1,
\label{minima3}
\end{equation}
where
\begin{equation}
\Delta_{\nu}(k_{min})=\frac{J_{\nu}(k_{min}\sigma)}  {\sigma^{\nu}} -
k_{min}\frac{J_{\nu-1}(k_{min})}{d}.
\label{minima4}
\end{equation}
The solution of equation (\ref{minima3}) produces the desired optimal values
of $\sigma_*$ and $\phi_*$, where
\begin{equation}
\phi_*=\frac{k_{min}^{\nu}}{ 2^{3\nu}\Gamma(1+\nu) \sigma_*^{2\nu} \Delta_{\nu}(k_{min})}.
\label{opt-phi}
\end{equation}

We find the solutions of (\ref{minima2}) by linearizing each Bessel function in (\ref{minima2})
around its respective first positive zero, i.e.,
\begin{eqnarray}
J_{\nu}(x)& =&\beta_1(x_0)(x-x_0)+{\cal O}\left((x-x_0)^2\right), \\
J_{\nu+1}(x)&=&\beta_2(y_0)(x-y_0)+{\cal O}\left((x-y_0)^2\right),
\end{eqnarray}
where
\begin{eqnarray}
\beta_1(x_0)&=&\frac{1}{2}[J_{\nu-1}(x_0)-J_{\nu+1}(x_0)], \label{beta1}\\
\beta_2(y_0)&=&\frac{1}{2}[J_{\nu}(y_0)-J_{\nu+2}(y_0)],
\end{eqnarray}
and $x_0$ and $y_0$ denote
the locations of the first positive zeros of $J_{\nu}(z)$ and $J_{\nu+1}(z)$, respectively.
Similarly, we employ the linearized form
\begin{equation}
xJ_{\nu}(x) =x_0\beta_1(x_0)(x-x_0)+{\cal O}\left((x-x_0)^2\right).
\end{equation}
Use of these relations in (\ref{minima2}) yields
the following equation for $k_{min}$:
\begin{equation}
k_{min} \simeq x_0-\frac{d(y_0-\sigma x_0)}{\frac{\beta_1}{\beta_2}\sigma^{\nu-1}x_0-d\sigma}.
\label{kmin}
\end{equation}
This formula provides an excellent approximation for $k_{min}$. For example, 
for $d=200$ (or $\nu=100$), substitution of the exact
values $x_0=108.8361659$, $y_0=109.8640469$ and $\beta_1/\beta_2=1.003189733$
as well as the numerical search solution $\sigma_*=1.008510$
into this formula predicts $k_{min}=108.4368917$.
This value is to be compared to the
numerical search solution of $k_{min}=108.4395$. This supports the fact that
the higher-order terms in the aforementioned linearized forms of the Bessel functions are negligibly small.
Indeed, we expect that this can be rigorously proved, but we shall not do so here. We will
assume the validity of the linearized forms in the asymptotics displayed below.

For large $d=2\nu$, we make use of the asymptotic formulas
\begin{eqnarray}
x_0&=&\nu+a_1\nu^{1/3}+\frac{a_2}{\nu^{1/3}}+\frac{a_3}{\nu} + {\cal O}(\frac{1}{\nu^{5/3}}), \\
y_0&=&\nu+a_1\nu^{1/3}+1+\frac{a_2}{\nu^{1/3}}+
\frac{a_2}{3\nu^{2/3}}+\frac{a_3}{\nu} + {\cal O}(\frac{1}{\nu^{4/3}}),
\end{eqnarray}
where the constants $a_1$, $a_2$ and $a_3$ are explicitly given in the Appendix.
For  $d=200$, these formulas predict $x_0=108.8362067$ and $y_0=109.8640871$, which
are in excellent agreement with the exact values reported in the preceding paragraph.
Using the asymptotic results given in  the Appendix, we obtain that
\begin{equation}
\frac{\beta_1}{\beta_2}= 1+\frac{2}{3\nu}-\frac{2C_2}{3 C_1\nu^{5/3}}+
{\cal O}\left(\frac{1}{\nu^2}\right),
%1+\frac{2}{3\nu}+\frac{0.981698886}{\nu^{5/3}}+{\cal O}\left(\frac{1}{\nu^2}\right).
\end{equation}
where the constants $C_1$ and $C_2$ are given explicitly
%$C_1=-1.104938082$ and $C_2=1.627074727$ 
in terms of the constants $a_1$ and $a_2$ in the Appendix.
For $d=200$, for example, this formula together with (\ref{C}) provide the estimate
$\beta_1/\beta_2 = 1.007122331$, which is to be compared to the exact 
result $\beta_1/\beta_2 = 1.006215695$.

The optimized asymptotic form for $\sigma_*$ is obtained by taking
the derivative of both sides of the zero-condition (\ref{minima3}) with respect to $\sigma$
and solving for $\sigma$ using relation (\ref{kmin}) for $k_{min}$. We obtain that
\begin{equation}
\sigma_*=1+\frac{q_1}{\nu}+\frac{q_2}{\nu^{5/3}}+{\cal O}\left(\frac{1}{\nu^2}\right),
\label{sigma*}
\end{equation}
where 
\begin{equation}
q_1=0.90763589355\ldots
\label{q1}
\end{equation}
is the unique positive root of $x e^x+ e^{2x}-5e^x+4=0$ and 
\begin{equation}
q_2=  \frac{a_1(8e^{q_1}-2q_1e^{q_1}-10e^{2q_1}+4+e^{3q_1}+4q_1 e^{2q_1})}
{3e^{q_1}(2q_1e^{q_1}-2q_1+12+3e^{2q_1}-13e^{q_1})}=-1.279349474.
\label{q2}
\end{equation}
Therefore, expression (\ref{kmin}) for $k_{min}$ has the asymptotic form
\begin{equation}
k_{min}=\nu+a_1\nu^{1/3}+Q_1+\frac{a_2}{\nu^{1/3}} +
{\cal O}\left(\frac{1}{\nu^{2/3}}\right),
\label{k}
\end{equation}
where
\begin{equation}
Q_1=\frac{2(q_1-1)}{e^{q_1}-2}=-0.3860921576.
\label{Q1}
\end{equation}
These formulas predict $\sigma_*=1.008482538$ and $k_{min}=108.4501542$, which
again are in excellent agreement with values reported above.

Linearizing each  Bessel function appearing in (\ref{minima4}) about
its first positive zero and using the results of the Appendix yields
\begin{eqnarray}
\Delta_{\nu}(k_{min})\simeq \frac{\beta_1(x_0)(k_{min}\sigma_*-x_0)}{\sigma_*^\nu}-
\frac{\beta_3(z_0)k_{min}(k_{min}-z_0)}{2\nu},
\end{eqnarray}
where $\beta_1(x_0)$ is given by (\ref{beta1}), $\beta_3(z_0)=[J_{\nu-2}(z_0)-J_{\nu}(z_0)]/2$,
and $z_0$ is the first positive zero of $J_{nu-1}$. Using relations (\ref{sigma*}) and (\ref{k}),
and the results of the Appendix yields the asymptotic expansion of $\Delta_{\nu}(k_{min})$:
\begin{eqnarray}
\Delta_{\nu}(k_{min})&=& \frac{D_1}{\nu^{2/3}} +
\frac{D_2}{\nu^{4/3}}+{\cal O}\left(\frac{1}{\nu^{5/3}}\right)
\end{eqnarray}
where
\begin{equation}
D_1=\frac{C_1(2-e^{q_1})}{2e^{q_1}}, 
\quad D_2=\frac{C_1\left[a_1(2e^{q_1}+6q_1e^{-q_1}-7)+3q_2(q_1-1)\right]}{3(2-e^{q_1})}+\frac{C_2 D_1}{C_1}.
\label{D}
\end{equation}
%Substitution of the constants $C_1,C_2,a_1,q_1$ and $q_2$ into the expressions above
%yields the numerical values $D_1=0.106651984$ and $D_2=0.336081209$.
%\begin{eqnarray}
%\Delta_{\nu}(k_{min})&=& \frac{0.106651984}{\nu^{2/3}} +
%\frac{0.336081209}{\nu^{4/3}}+{\cal O}\left(\frac{1}{\nu^{5/3}}\right)
%\end{eqnarray}
Note also that
\begin{eqnarray}
\left(\frac{k_{min}}{\nu}\right)^{\nu} &=& e^{a_1\nu^{1/3}+Q_1}
\Bigg[1+\frac{E_1}{\nu^{1/3}}+
\frac{E_2}{\nu^{2/3}}+{\cal O}\left(\frac{1}{\nu}\right)\Bigg] \nonumber \\
%1-\frac{0.688767207}{2\nu^{1/3}}
%+\frac{0.953693390}{\nu^{2/3}}+{\cal O}\left(\frac{1}{\nu}\right)\Bigg] \nonumber \\
%&=&e^{1.8557571\nu^{1/3}-0.3860921576}
%\Bigg[1-\frac{0.688767207}{\nu^{1/3}}
%+\frac{0.9536932716}{\nu^{2/3}}+{\cal O}\left(\frac{1}{\nu}\right)\Bigg]
\end{eqnarray}
and
\begin{eqnarray}
\sigma_*^{2\nu} &=& e^{2q_1}\left[1+\frac{2q_2}{\nu^{2/3}}-\frac{q_1^2}{\nu}+\frac{2q_2^2}{\nu^{4/3}}+
{\cal O}\left(\frac{1}{\nu^{5/3}}\right)\right],
%&=&6.142745469-\frac{15.71743637}{\nu^{2/3}}-\frac{5.060411626}{\nu}
%+\frac{20.10809396}{\nu^{4/3}}+{\cal O}\left(\frac{1}{\nu^{5/3}}\right)
\end{eqnarray}
where 
\begin{equation}
E_1=a_2-\frac{a_1^2}{2}, \quad E_2=-Q_1a_1+\frac{a_1^4-4a_1^2a_2+4a_2^2}{8},
\label{E}
\end{equation}
and $Q_1$ is given by (\ref{Q1}). For $d=200$, these formulas [together
with the constants specified by (\ref{q1}), (\ref{q2}), (\ref{D}), \ref{E}), (\ref{a}) and (\ref{C})] 
predict $\Delta_{\nu}(k_{min})=0.00567441932$, 
$(k_{min}/\nu)^\nu=3353.018128$ and $\sigma_*^{2\nu}=5.405924156$. These
values should be compared to the exact value of
$\Delta_{\nu}(k_{min})=0.00559813885$, $(k_{min}/\nu)^\nu=3301.799093$
and $\sigma_*^{2\nu}=5.445550297$.

Thus, substituting the asymptotic relations above into the optimal expression (\ref{opt-phi}) for the
density and invoking Conjecture \ref{conj} yields the conjectural lower bound}
\begin{eqnarray}
\phi_{\mbox{\scriptsize max}} \ge \phi_*&=& 
 \frac{1}{2^{[3-\log_2(e)]\nu-\log_2(e)a_1\nu^{1/3}+(2q_1-Q_1)\log_2(e)}}
\Bigg[\frac{1}{2D_1}\sqrt{\frac{2}{\pi}}\Big[\nu^{1/6}+\frac{E_1}{\nu^{1/6}} \nonumber \\
&+& \frac{E_2-2q_2-D_2/D_1}{\nu^{1/2}}+ {\cal O}\left(\frac{1}{\nu^{5/6}}\right)\Big]\Bigg], 
%&=& 
% \frac{1}{2^{1.557304959\nu-2.677291565\nu^{1/3}+0.5570132411}}
%\Bigg[0.6089457595 \nu^{1/6}- \frac{0.4194218702}{\nu^{1/6}} \nonumber \\&&
%+\frac{0.06159502415}{\nu^{1/2}}
%+ {\cal O}\left(\frac{1}{\nu^{5/6}}\right)\Bigg], 
\label{term3}
\end{eqnarray}
where we have used the asymptotic relation $\Gamma(1+\nu)\sim\nu^{\nu}\sqrt{2\pi\nu} e^{-\nu}$.
For $d=200$, this asymptotic formula [together
with the constants specified by (\ref{q1}), (\ref{q2}), (\ref{D}), \ref{E}), (\ref{a}) and (\ref{C})] 
predicts $\phi_*=5.626727001  \times 10^{-44}$, which
is in good agreement with the numerical search solution of $\phi_*=5.667098 \times 10^{-44}$.
Note also that the formula (\ref{opt-phi}) with $k_{min}$ estimated from (\ref{kmin}) yields
$\phi_*=5.666392126 \times 10^{-44}$, which is remarkably close to the numerical
solution. For large $d$,  result (\ref{term3}) yields the following dominant
asymptotic formula for the conjectural lower bound on $\phi_{\mbox{\scriptsize max}}$:
\begin{equation}
\phi_{\mbox{\scriptsize max}} \ge 
\phi_* \sim \frac{d^{1/6}}{2^{2/3}D_1\sqrt{\pi}\,2^{[3-\log_2(e)]d/2}}=
\frac{3.276100896 ~d^{1/6}}{2^{0.7786524795 \ldots \,d}}.
\label{nice-result}
\end{equation}
This putatively provides the long-sought exponential improvement on Minkowski's
lower bound. Note that the constant $D_1=0.1084878572$  appearing
in (\ref{nice-result}) is determined
from the appropriate relation in  (\ref{D}) using the value for $q_1$ given by (\ref{q1}) and the more
refined estimate of $C_1$ given by (\ref{C11}). 

Substitution of the asymptotic expression (\ref{nice-result}) into (\ref{Z}) and use of (\ref{true-bound-Z}) yields a
conjectural lower bound on the maximal kissing number
\begin{equation}
Z_{\mbox{\scriptsize max}} \ge Z_* \sim \frac{2^{1/3}e^{2q_1}d^{1/6}}{D_1\sqrt{\pi}}\, 2^{[\log_2(e)-1]d/2}=
40.24850787 \,d^{1/6} \, 2^{0.2213475205\ldots\, d},
\end{equation}
which applies for large $d$.
This result is superior to the best known 
asymptotic lower bound on the maximal kissing number of $2^{0.2075\ldots d}$
\cite{Wy65}. Note that such a disordered packing would be substantially {\it hyperstatic} (the average kissing
number is greater than $2d$ \cite{Do05c}) and therefore would be appreciably
different from a maximally random jammed packing \cite{To00b,To02a}, which is isostatic (see Section 4.3)
and hence signficantly smaller in density.

\section{Discussion}
\label{discuss}

Our results have immediate implications for the linear programming bounds
of Cohn and Elkies \cite{Co03}, regardless of
the validity of Conjecture \ref{conj}.  As we noted earlier, our optimization
procedure is precisely the dual of the their primal linear programming
upper bound. The existence of our test functions (\ref{step2})
and (\ref{step3}) that satisfy the conditions (i), (ii) and (iii) of Section 5 for 
densities up to the terminal density $\phi_*$ narrows the duality gap [cf. (\ref{term2})
and (\ref{term3})]. In particular, inequality (\ref{term3}) provides the greatest lower bound
known for the dual linear program. Moreover, the existence of the inequalities (\ref{lower-2})
and (\ref{term3}) proves that linear programming bounds cannot possibly match the Minkowski 
lower bound for any dimension $d$. Finally, this link to the Cohn-Elkies formulation
proves that the terminal density $\phi_*$ can never exceed the Cohn-Elkies upper bound, 
which obviously must be true if the terminal density corresponds to a rigorous lower bound.

Conjecture \ref{conj} concerning the existence  of disordered sphere packings 
is plausible for a number of reasons: (i) the decorrelation principle of Section
4.3; (ii) the necessary Yamada condition appears to only have relevance in very low dimensions;
(iii) other new necessary conditions described in Section 5.4
also seem to be germane only in very low dimensions; (iv) 
there are numerical constructions of configurations of disordered sphere
packings on the torus corresponding to these test radial distributions functions
in low dimensions for densities up to the terminal density \cite{Cr03,Ou06}; and 
(v) the test radial distributions functions (\ref{step1}) and (\ref{step2})
recover the asymptotic forms of known rigorous bounds. Concerning the latter point,
if Conjecture \ref{conj} is false, it is certainly not revealed by the results
produced by the test functions (\ref{step1}) and (\ref{step2}) because the forms of known
rigorous results, obtained using completely different techniques, are recovered.
If Conjecture \ref{conj} is false, why would it suddenly be revealed by the introduction of a gap in the 
test radial distribution function [cf. (\ref{step3})] relative to (\ref{step2})?
This would seem to be unlikely and lends credibility to the conjecture in our view.

Conjecture \ref{conj}, the particular choice (\ref{step3}) and our optimization procedure leads to 
a lower bound on the maximal density 
that improves on the Minkowski bound by an exponential factor. 
Our results suggest that the densest packings in sufficiently
high dimensions may be disordered rather than periodic, implying
the existence of disordered classical ground states for some continuous 
potentials. A byproduct of our 
procedure is an associated lower bound on the maximal kissing number, which is superior 
to the currently best known result. By no means is the choice (\ref{step3}) optimal. 
For example, one may be able to improve our putative lower bound by allowing the 
test function to be some positive function smaller than unity for $1 \le r \le \sigma$. Of course,
it would be desirable to choose test functions that make asymptotic
analysis exactly rather than numerically tractable.

Our putative exponential improvement over Minkowksi's bound
is striking and  should provide some impetus to determine
the validity of Conjecture \ref{conj}. As a first step in this
direction, it would be fruitful if one could  show that
for sufficiently small densities, the two standard nonnegativity conditions
on the pair correlation function $g_2$ are sufficient to
assure the existence of a point process, whether it is a sphere packing
or not. Another problem worth pursuing is the demonstration of
the existence of a pair interaction
potential in $\mathbb{R}^d$ corresponding to a sphere packing for a given $\rho$ and $g_2$
provided that $\rho$ and $g_2$ are sufficiently small. Such a proof may be possible
following the methods that Koralov \cite{Ko05} used for the
lattice setting. It would also be profitable to pursue the construction of disordered sphere
packings with densities that exceed $1/2^d$ for sufficiently large $d$.

\noindent{\bf Acknowledgments}
\smallskip

\noindent This work supported in part by the National Science Foundation
under Grant No. DMS-0312067.
We thank Peter Sarnak and Enrico Bombieri for encouraging us
to apply our methods to try to obtain exponential improvement
on Minkowski's lower bound. We benefited greatly from discussions
with Henry Cohn, John Conway, Leonid Koralov and Thomas Spencer.

\section*{APPENDIX}
\renewcommand{\theequation}{A-\arabic{equation}}
 \setcounter{equation}{0}  % reset counter 

An asymptotic expression for $J_{\nu}(x)$ when $\nu$ is large and $x >\nu$ is given by \cite{Wa58}
\begin{equation}
J_{\nu}(x) = A_{\nu}(x) \left[\cos[\omega_{\nu}(x)-\pi/4]+ {\cal O}\left( \frac{3x^2+2\nu^2}{12(x^2-\nu^2)}\right) \right],
\label{watson}
\end{equation}
where
\begin{equation}
A_{\nu}(x)=\Bigg[\frac{2}{\pi \sqrt{x^2-\nu^2}}\Bigg]^{1/2}
\end{equation}
and 
\begin{equation}
\omega_{\nu}(x)=\sqrt{x^2-\nu^2}-\nu \cos^{-1}(\nu/x).
\end{equation} 
The function $A_{\nu}(x) \cos[\omega_{\nu}(x)-\pi/4]$  in (\ref{watson}) actually represents the dominant term
in the asymptotic expansion (4) on page 244 of Watson \cite{Wa58} for $J_{\nu}(x)$ when $\nu$ is large and $x >\nu$
and $A_{\nu}(x)$ multiplied by the big-$\cal O$ term  represents the largest absolute error when this dominant 
term is used to estimate $J_{\nu}(x)$. A problem of central concern is an estimate of $J_{\nu}(x)$ in the vicinity of its
first positive zero $x_0$ when $\nu$ is large. The first positive zero 
has the asymptotic expansion \cite{Ol60}
\begin{equation}
x_0=\nu+a_1\nu^{1/3}+\frac{a_2}{\nu^{1/3}}+\frac{a_3}{\nu} + {\cal O}(\frac{1}{\nu^{5/3}}),
\label{zero}
\end{equation}
where 
\begin{equation}
a_1=1.8557571\ldots, \qquad a_2= 1.033150\ldots, \qquad \;  a_3=-0.003971\ldots. 
\label{a}
\end{equation}
Expanding $J_{\nu}(x)$ in a Taylor series about $x=x_0$ and neglecting quadratic
and higher-order terms gives the linear estimate
\begin{equation}
J_{\nu}(x)\simeq \frac{1}{2}[J_{\nu-1}(x_0)-J_{\nu+1}(x_0)](x-x_0),
\end{equation}
where we take the Bessel functions on the right side to be given by 
the asymptotic forms
\begin{equation}
J_{\nu+1}(x_0) = A_{\nu+1}(x_0) \left[\cos[\omega_{\nu+1}(x_0)-\pi/4]+ {\cal O}\left( \frac{3x^2+2\nu^2}{12(x^2-\nu^2)}\right) \right]
\end{equation}
and 
\begin{equation}
J_{\nu-1}(x_0) = A_{\nu-1}(x_0) \left[\cos[\omega_{\nu-1}(x_0)-\pi/4]+ {\cal O}\left( \frac{3x^2+2\nu^2}{12(x^2-\nu^2)}\right)\right].
\end{equation}
We will also need to consider the related functions
\begin{equation}
J_{\nu+1}(x)\simeq\frac{1}{2}[J_{\nu}(y_0)-J_{\nu+2}(y_0)](x-y_0),
\end{equation}
\begin{equation}
J_{\nu-1}(x)\simeq\frac{1}{2}[J_{\nu-2}(z_0)-J_{\nu}(z_0)](x-z_0),
\end{equation}
where $y_0$ and $z_0$ are the first positive zeros of $J_{\nu+1}(x)$ and
$J_{\nu-1}(x)$, respectively, which are asymptotically given by
\begin{eqnarray}
y_0&=&\nu+a_1\nu^{1/3}+1+\frac{a_2}{\nu^{1/3}}+\frac{a_2}{3\nu^{2/3}}+\frac{a_3}{\nu} + {\cal O}(\frac{1}{\nu^{4/3}}),\\
z_0&=&\nu+a_1\nu^{1/3}-1+\frac{a_2}{\nu^{1/3}}-\frac{a_2}{3\nu^{2/3}}+\frac{a_3}{\nu} + {\cal O}(\frac{1}{\nu^{4/3}}).
\end{eqnarray}
Note that the asymptotic expressions for the zeros given here for $\nu=100$ ($d=200$) predict $x_0=108.8362071$,
$y_0=109.8641774$ and $z_0=107.8082369$, which are in excellent agreement with
exact results $x_0=108.8361659$, $y_0=109.8640469$ and $z_0=107.8081033$ as obtained from Maple.

Using Maple and the results above, we obtain the following asymptotic expansions:
\begin{eqnarray}
\frac{1}{2}[J_{\nu-1}(x_0)-J_{\nu+1}(x_0)]&=&\frac{C_1}{\nu^{2/3}}+\frac{C_2}{\nu^{4/3}}
+{\cal O}\left(\frac{1}{\nu^{2}}\right), \\
\frac{1}{2}[J_{\nu}(y_0)-J_{\nu+2}(y_0)]&=&
\frac{C_1}{\nu^{2/3}}+\frac{C_2}{\nu^{4/3}}
-\frac{2C_1}{3\nu^{5/3}}+{\cal O}\left(\frac{1}{\nu^{2}}\right), \\
\frac{1}{2}[J_{\nu-2}(z_0)-J_{\nu}(z_0)] &=&\frac{C_1}{\nu^{2/3}}+\frac{C_2}{\nu^{4/3}}
+\frac{2C_1}{3\nu^{5/3}}+{\cal O}\left(\frac{1}{\nu^{2}}\right).
\end{eqnarray}
where
\begin{eqnarray}
\hspace{-0.6in} C_1&=& -\frac{2^{1/4}\left[\sqrt{2}f_1(a_1)+8a_1^{3/2}f_2(a_1)\right]}{8\sqrt{\pi}a_1^{5/4}}, \label{C1-1}\\
\hspace{-0.6in} C_2&=&\frac{2^{3/4}\left[1152  a_1^6 -3840  a_1^4 a_2-180 a_1^3 +600a_1a_2 -225
\right]f_1(a_1)+2^{1/4}\left[3072 a_1^{9/2}-200a_1^{3/2}\right]f_2(a_1)}{3840\sqrt{\pi}a_1^{13/4}}, \nonumber\\ \\
%C_3&=&\frac{2^{3/4}\left[-288 a_1^{15/2}+384 a_1^{11/2}a_2+1920  a_1^{7/2}a_2^2 +150 a_1^{3/2}\right]f_1(a_1)}
%{7680\sqrt{\pi}a_1^{15/4}} \nonumber \\
%&-& \frac{2^{1/4}\left[360a_1^6+1440(a_1^4a_2+a_1^2a_2^2)-290a_1^3+1920\sqrt{2}a_1^{5/2}a_2+380a_1a_2-275\right]
%f_2(a_1)}{7680\sqrt{\pi}a_1^{15/4}},\\
f_1(a_1)&=&\sin\left(\frac{(2a_1)^{3/2}}{3}\right)+\cos\left(\frac{(2a_1)^{3/2}}{3}\right), \qquad
f_2(a_1)=\sin\left(\frac{(2a_1)^{3/2}}{3}\right)-\cos\left(\frac{(2a_1)^{3/2}}{3}\right).
\end{eqnarray}
Thus, substitution of the values for the constants $a_1$ and $a_2$ into
the expressions above yields the estimates 
\begin{equation}
C_1=-1.104938082, \qquad C_2=1.627074727.
\label{C}
\end{equation}
For $d=200$ ($\nu=100$), for example, the asymptotic expansions (A-13) - (A-15) predict 
$-0.04778125640$,  $-0.04743934518$ and $-0.04812316762$, respectively.
This is to compared to the corresponding exact results: 
$-0.04829366129$,   $-0.04799533693$ and $-0.04859672879$.
Note that although the estimates of the constants $C_1$ and $C_2$ given by (\ref{C})
involve a small error due to our use of only the dominant asymptotic term (\ref{watson}),
the functional forms of the asymptotic expansions (A-13) - (A-15) are exact.

We can show that the exact expressions for the constants $C_1$ and $C_2$ appearing in (A-13) - (A-15)
are rapidly converging asymptotic expansions in inverse powers of the 
constant $a_1$ appearing in (\ref{zero}).
For example, using expansion (4) on page 244 of Watson, we find that the first three
terms of the expansion of $C_1$ is given by
\begin{equation}
C_1=\frac{C_{11}}{a_1^{5/4}}+\frac{C_{12}}{a_1^{11/4}}+\frac{C_{13}}{a_1^{17/4}} +{\cal O}\left( \frac{1}{a_1^{23/4}}\right),
\label{C1-2}
\end{equation}
where
\begin{displaymath}
C_{11}=-\frac{2^{1/4}\left[\sqrt{2}f_1(a_1)+8a_1^{3/2}f_2(a_1)\right]}{8\sqrt{\pi}}, \quad
C_{12}=\frac{5 \cdot 2^{1/4}\left[4\sqrt{2}f_1(a_1)-7a_1^{3/2}f_2(a_1)\right]}{384\sqrt{\pi}},
\end{displaymath}
\begin{displaymath}
C_{13}=-\frac{385\cdot 2^{1/4}\left[13\sqrt{2}f_1(a_1)+8a_1^{3/2}f_2(a_1)\right]}{221184\sqrt{\pi}}.
\end{displaymath}
The first term of the expansion (\ref{C1-2}) is the dominant one and is identical to the estimate given
in (\ref{C1-1}). The first term of (\ref{C1-2}) is about 66 times larger than the second
term in absolute value and the second term is about 7 times larger than the third term in absolute value. 
Substitution of the constant $a_1$ into (\ref{C1-2}) yields the more refined estimate
\begin{equation}
C_1=-1.123958144.
\label{C11}
\end{equation}
This refined estimate differs from the dominant first-term estimate given
in (\ref{C}) in the third significant figure.
One could continue correcting this estimate by including additional terms
in the asymptotic expansion but this quickly  becomes tedious and is not necessary because,
as we show at the end of Section 5, the precise value of $C_1$ is not relevant for
the putative exponential improvement of Minkowski's lower bound on the density,
as specified by (\ref{nice-result}).

%\begin{eqnarray}
%\frac{1}{2}[J_{\nu-1}(x_0)-J_{\nu+1}(x_0)]&=&-\frac{1.104938082}{\nu^{2/3}}+\frac{1.627074726}{\nu^{4/3}}
%-\frac{1.929036634}{\nu^2}+{\cal O}\left(\frac{1}{\nu^{8/3}}\right), \\
%\frac{1}{2}[J_{\nu}(y_0)-J_{\nu+2}(y_0)]&=&-\frac{1.104938082}{\nu^{2/3}}+\frac{1.627074726}{\nu^{4/3}}
%+\frac{0.7366253880}{\nu^{5/3}}\nonumber\\&&
%-\frac{1.929036634}{\nu^2}+{\cal O}\left(\frac{1}{\nu^{7/3}}\right),\\
%\frac{1}{2}[J_{\nu-2}(z_0)-J_{\nu}(z_0)]&=&-\frac{1.104938082}{\nu^{2/3}}+\frac{1.627074726}{\nu^{4/3}}
%-\frac{0.7366253880}{\nu^{5/3}}.\nonumber \\&&
%-\frac{1.929036634}{\nu^2}+{\cal O}\left(\frac{1}{\nu^{7/3}}\right).
%\end{eqnarray}

\vspace{-0.2in}

%\bibliographystyle{alpha}
%\bibliography{new}

\end{document}